\documentclass[11pt]{article}
\usepackage{makeidx}
\usepackage{amsmath,amssymb}
\usepackage{verbatim}

\usepackage{exscale}
\usepackage{tabularx}
\usepackage{amsthm}
\usepackage{latexsym}
\usepackage[usenames,dvipsnames]{color}
\usepackage[dvips]{graphicx,epsfig}
\usepackage{hyperref}
\usepackage{url}

\newcommand{\mb}[1]{\mbox{\boldmath $#1$}}

\newenvironment{pf}{\noindent {\bf Proof. }}{\hfill $\square$}

\newcommand{\bea}{\begin{eqnarray*}}
\newcommand{\eea}{\end{eqnarray*}}

\newcommand{\req}[1]{$(\ref{#1})$}

\newcommand{\vect}[1]{{\boldsymbol #1 }}

\newcommand{\inprod}[2]{\langle #1 , #2 \rangle }
\newcommand{\card}[1]{\left\lvert #1 \right\rvert}
\newcommand{\bc}{\begin{center}}
\newcommand{\ec}{\end{center}}

\newcommand{\by}{\vect{y}}
\newcommand{\bY}{\vect{Y}}
\newcommand{\refs}[1]{$(\ref{#1})$}

\newcommand{\bX}{\vect{X}}
\newcommand{\bx}{\vect{x}}

\newcommand{\bZ}{\boldsymbol Z}
\newcommand{\R}{\mathbb R}

\newcommand{\be}{\begin{equation}}
\newcommand{\ee}{\end{equation}}
\newcommand{\beaa}{\begin{eqnarray*}}
\newcommand{\eeaa}{\end{eqnarray*}}
\newcommand{\ben}{\begin{enumerate}}
\newcommand{\een}{\end{enumerate}}
\newcommand{\db}{\hspace*{\fill}{\zapf o}}

\newcommand{\bpn}{\begin{proposition}\twlsf}
\newcommand{\epn}{\db\end{proposition}}
\newcommand{\bdm}{\begin{displaymath}}
\newcommand{\edm}{\end{displaymath}}
\newcommand{\ba}{\begin{array}}
\newcommand{\ea}{\end{array}}

\newcommand{\st}{\mathop{\rm s.t.}}

\newtheorem{lemma}{Lemma}
\newtheorem{proposition}{Proposition}
\newtheorem{corollary}{Corollary}

\newtheorem{theorem}{Theorem}

\newcommand{\eps}{\epsilon}

\newcommand{\norm}[1]{\left\lVert#1\right\rVert}
\newcommand{\tnorm}[1]{\lVert\mkern-2mu |#1|\mkern-2mu\rVert}
\title{Finding approximately rank-one submatrices with the nuclear norm
  and $\ell_1$-norm\thanks{Supported in part by the U.~S.~Air Force 
Office of Scientific Research, a Discovery Grant from the 
Natural
Sciences and Engineering Research Council of Canada, and a grant
from MITACS.}}
\author{Xuan Vinh Doan\thanks{Department of
    Combinatorics and Optimization, University of Waterloo, 200
    University Avenue West, Waterloo, ON N2L 3G1, Canada,
    vanxuan@uwaterloo.ca.} \and Stephen Vavasis\thanks{Department of
    Combinatorics and Optimization, University of Waterloo, 200
    University Avenue West, Waterloo, ON N2L 3G1, Canada,
    vavasis@uwaterloo.ca.}}

\date{November 2010}

\begin{document}
\maketitle

\begin{abstract}
We propose a convex optimization formulation with the nuclear norm and
$\ell_1$-norm to find a large approximately rank-one submatrix of a given
nonnegative matrix. We
develop optimality conditions for the formulation and characterize the
properties of the optimal solutions. We establish conditions under
which the optimal solution of the convex
formulation has a specific sparse
structure. Finally, we show that, under
certain hypotheses, with high probability,
the approach can recover the rank-one submatrix even when it is 
corrupted with random noise and inserted as a submatrix into a much
larger random noise matrix.
\end{abstract}

\section{Introduction}
\label{sec:intro}
Given a nonnegative matrix $\mb{A}\in\R^{m\times n}$, $a_{ij}\geq 0$
for all $i=1,\ldots,m$, $j=1,\ldots,n$, we consider the problem
of finding $\mathcal{I}\subset\{1,\ldots,m\}$ and
$\mathcal{J}\subset\{1,\ldots,n\}$ such that 
$\mb{A}(\mathcal{I},\mathcal{J})$ is close to a rank-one matrix,
and such that 
$\Vert \mb{A}(\mathcal{I},\mathcal{J})\Vert$ is large.
We shall call this problem the LAROS problem
(for ``large approximately rank-one submatrix'').

The main application of the LAROS problem is for finding features in data.  
For example, suppose $\mb{A}$ represents a corpus of documents
in some language.
Each column of $\mb{A}$ is in correspondence with one document,
and each
row is in correspondence with a term used in the corpus.  Here,
``term'' means a word in the
language, excluding common words such
as articles and prepositions.  The $(i,j)$ entry of $\mb{A}$ is the 
number of occurrences
of term $i$ in document $j$, perhaps normalized.  Such a matrix
is called the {\em term-document matrix} of the underlying corpus.

In this case, an approximately rank-one submatrix of $\mb{A}$
corresponds to a subset of terms and a subset of documents in which
the selected terms occur with proportional frequencies in the selected
documents.  Such a submatrix may correspond to the intuitive notion
of a topic that recurs in several documents, since a topic may manifest
itself as a particular group of relevant terms that occur roughly in
the same proportions.

As another example, the matrix $\mb{A}$ may correspond to a
database of pixelated grayscale images, where each image has
the same pixel size. Each column of $\mb{A}$
corresponds to one image, and each
row to one pixel position.  The $(i,j)$ entry of $\mb{A}$ is the
intensity of the $i$th pixel in the $j$th image.  In this case, the
approximately 
rank-one submatrix corresponds to a visual
feature that recurs in a certain
position in some subset of the images.

If one wanted to find more than one topic in a term-document matrix
or more than one feature in an image database matrix, then one
could iteratively find an approximately rank-one submatrix, subtract
it from $\mb{A}$ (perhaps modifying the result of the subtraction
to ensure that $\mb{A}$ remains nonnegative), 
and then repeat the procedure $p$
times.  Let the submatrices discovered
be denoted $(\mathcal{I}_1,\mathcal{J}_1)$, \ldots,
$(\mathcal{I}_p,\mathcal{J}_p)$.  Suppose
$A(\mathcal{I}_i,\mathcal{J}_i)\approx \vect{w}_i\vect{h}_i^T$ for
$i=1,\ldots, p$, and let $\bar{\vect{w}}_i$,
$\bar{\vect{h}}_i$ denote the extension of
$\vect{w}_i,\vect{h}_i$ to vectors of
length $m,n$ by inserting zeros for entries
not in 
$\mathcal{I}_i,\mathcal{J}_i$ respectively.

It is known that
if $\hat{\mb{A}}$ is a nonnegative matrix representing a submatrix
of $\mb{A}$, then
the minimizer $\vect{w}\vect{h}^T$ of
$\Vert \hat{\mb{A}}-\vect{w}\vect{h}^T\Vert$ is the dominant singular
vector pair in either the Frobenius
or 2-norm (a consequence of the Eckart--Young theorem, Theorem
2.5.3 of \cite{GVL}) and furthermore, $\vect{w}\ge\vect{0}$ and
$\vect{h}\ge\vect{0}$ (a consequence of the
Perron--Frobenius theorem.)
Thus, without loss of generality, we may assume that each
$\vect{w}_i\vect{h}_i^T$ determined by the iterative
computation is nonnegative.

In this case, one has
an approximate factorization 
$$\mb{A}\approx[\bar{\vect{w}}_1,\ldots,
\bar{\vect{w}}_p][\bar{\vect{h}}_1,\ldots,\bar{\vect{h}}_p]^T,$$
where we can write the right-hand side as $\mb{W}\mb{H}^T$ with
$\mb{W}\ge\mb{0}$, $\mb{H}\ge\mb{0}$.
This factorization
is called a {\em nonnegative matrix factorization} of $\mb{A}$.
The earliest reference known to us concerning nonnegative
matrix factorization is Thomas' solution \cite{Thomas} to a problem posed
by A.\ Berman and R.\ Plemmons (which, according to a remark in the
journal, was also solved by A.\ Ben-Israel).
Cohen and Rothblum \cite{CohenRothblum} describe applications
for NMF in probability, quantum mechanics and other fields.
Lee and Seung \cite{LeeSeung} showed that NMF can find features
in image databases, and Hofmann \cite{hofmann} showed that probabilistic
latent semantic analysis, a variant of NMF, can effectively
cluster documents according to their topics.  

Nonnegative matrix factorization is sometimes posed as an optimization
problem: find $\mb{W}\in\R^{m\times p}$ and $\mb{H}\in\R^{n\times p}$,
both nonnegative, such that $\Vert \mb{A}-\mb{W}\mb{H}^T\Vert$ is minimized in some matrix
norm.   It is known that this optimization problem is
NP-hard \cite{nmfnphard}.  Therefore, it is not surprising that
most algorithms for the problem are heuristic in the sense that they
do not make guarantees about the quality of the approximation.

One class of heuristic NMF algorithms are the `greedy' algorithms
\cite{greiner, bergman, r1d, gillisglineur} 
that follow the framework described above.
In a greedy
algorithm, the columns of $\mb{W}$ and $\mb{H}$ are generated sequentially,
with each new pair of columns accounting for one feature in the
original $\mb{A}$.  These greedy algorithms give rise to the LAROS
subproblem
addressed in this paper, namely, find one pair $\vect{w}_i,\vect{h}_i$
nonnegative such that $\vect{w}_i\vect{h}_i^T$, 
is a good approximation for a submatrix of $\mb{A}$ in the positions $(i,j)$
where $\mb{A}$ is positive.

The LAROS subproblem, however, is itself
NP-hard as observed by \cite{gillisglineur}.
This is because the maximum-edge
biclique problem can be naturally 
expressed as a rank-one submatrix problem. The biclique problem takes
as input a bipartite graph $G=(U,V,E)$.  The output is composed
of two subsets $U^*\subset U$ and $V^*\subset V$ such that
$U^*\times V^*\subset E$ (i.e., all possible $|U^*|\cdot|V^*|$ edges
between $U^*$ and $V^*$ are present in $G$) and such that
$|U^*|\cdot |V^*|$ is maximum with this property.  
This problem was shown by Peeters \cite{Peeters} to be NP-hard.

Maximum-edge biclique can be expressed as finding a large rank-one
submatrix using the following construction.  Let $\mb{A}$ be a
$|U|\times |V|$ matrix with rows in correspondence to $U$ and columns
in correspondence to $V$.  Entry $(i,j)$ of $\mb{A}$ for $(i,j)\in U\times
V$ is 1 if $(i,j)\in E$, else this entry is 0.  Then a biclique
corresponds exactly of a $|U^*|\cdot|V^*|$ submatrix of all 1's.  A
submatrix of all 1's is a rank-one matrix of norm
$(|U^*|\cdot|V^*|)^{1/2}$, and there is no other kind of rank-one
submatrix of $\mb{A}$.

We will provide a formal definition of the LAROS
problem, i.e., exactly what is the desired output in Section~\ref{sec:formulation}.  (Some
authors mentioned earlier, e.g., \cite{r1d} and \cite{gillisglineur} have
provided other formal definitions.)  We
will also propose an convex optimization problem in 
Section~\ref{sec:formulation} that,
for matrices $\mb{A}$ constructed in a certain way,
successfully finds large, approximately rank-one
submatrices.  
We present two such theorems.  One case is when the
approximately rank-one submatrix dominates the rest of the matrix;
this is presented in Section~\ref{sec:sparsity}.

The second case is when $\mb{A}$ is constructed as follows:
$$\mb{A}=\mb{A}_0+\mb{R},$$ 
where there are two index sets $\mathcal{I},\mathcal{J}$ such that
$\mathop{\rm rank}(\mb{A_0}(\mathcal{I},\mathcal{J}))=1$,
$\mb{A}_0(i,j)=0$ for all $(i,j)\notin \mathcal{I}\times \mathcal{J}$,
and $\mb{R}$ is a random matrix representing noise.  In this case,
under certain assumptions, our algorithm recovers
$(\mathcal{I},\mathcal{J})$ from $\mb{A}$ as proved in
Section~\ref{sec:randomnoise}.

\section{Matrix norm minimization}
\label{sec:formulation}

For reasons that will become clear, we start this section by presenting the
convex relaxation of the LAROS problem, and only later
will we present a nonconvex exact optimization formulation.
In particular, we propose the following convex optimization problem that
in some cases solves the LAROS problem:
\begin{equation}
\begin{array}{rl}
\min  & \Vert \mb{X}\Vert_*+\theta \Vert \mb{X}\Vert_1  \\
\mbox{s.t.} & \inprod{\mb{A}}{\mb{X}}\ge 1.
\end{array}
\label{eq:LAROSrelaxation}
\end{equation}
The matrix $\mb{X}\in\R^{m\times n}$ is the unknown.
Norm $\Vert\mb{X}\Vert_*$ is the {\em nuclear norm}, also called
the {\em trace norm}; it is the sum of the singular values of 
$\mb{X}$.  We use the notation $\Vert\mb{X}\Vert_1$ to mean the sum
of the absolute values of entries of $\mb{X}$, that is, the 
$\ell_1^{mn}$-norm applied to ${\rm vec}(\mb{X})$,  the concatenation
of the columns of $\mb{X}$ into a long vector.  Finally,
$\inprod{\mb{A}}{\mb{X}}$ means the inner product of the two matrices.
We note that an objective function involving a sum of the
nuclear and $\ell_1$-norms was used for a different
purpose by \cite{Chandrasekaran}.
Two other norms used extensively in this paper are
$\Vert\mb{X}\Vert$, which is the spectral or 2-norm, i.e.,
$\sigma_1(\mb{X})$, and $\Vert\mb{X}\Vert_\infty$, which is
the $\ell_\infty^{mn}$-norm applied to ${\rm vec}(\mb{X})$, i.e.,
the maximum absolute entry of $\mb{X}$.

Before beginning a detailed analysis of this 
optimization problem, we first provide
some motivation.  Consider first the simplification obtained by
taking $\theta=0$:
$$
\begin{array}{rl}
\min  & \Vert \mb{X}\Vert_*  \\
\mbox{s.t.} & \inprod{\mb{A}}{\mb{X}}\ge 1.
\end{array}
$$
It follows from Proposition~\ref{prop:nuclear} below that the optimal solution 
is found using the singular value decomposition.  
In particular, if $\mb{A}$ is factored as
$\mb{A}=\mb{U}\mb{\Sigma}\mb{V}^T$, where $\mb{U}\in\R^{m\times m}$ 
is orthogonal, $\mb{\Sigma}\in\R^{m\times n}$
is diagonal, and $\mb{V}\in\R^{n\times n}$ is orthogonal, then an optimizer is $\mb{X}=\mb{U}(:,1)\mb{V}(:,1)^T/\sigma_1$.
Thus, when $\theta=0$, the above formulation successfully
finds the best rank-one approximation to the whole matrix $A$.

This approximation, however, is not always well suited for identifying
submatrices.  Consider e.g. the following $6\times 6$ matrix:
$$
\mb{A}=
\left(
\begin{array}{cccccc}
0.8 & 0.9 & 1.1 & 0.1 & 0.2 & 0.2 \\
0.8 & 1.1 & 0.8 & 0 & 0 & 0 \\
1.0 & 1.0 & 0.8 & 0 & 0 & 0 \\
0 & 0 & 0 & 0.8 & 0.9 & 1.0 \\
0 & 0 & 0 & 0.9 & 1.0 & 0.8 \\
0 & 0 & 0 & 1.0 & 1.1 & 0.8
\end{array}
\right)
$$
It is apparent from inspection that this matrix has two $3\times 3$
approximately rank-one blocks in positions $\{1,2,3\}\times\{1,2,3\}$
and $\{4,5,6\}\times\{4,5,6\}$.  If the `noise' entries in the upper
right $\{1,2,3\}\times \{4,5,6\}$ block were absent, then the two dominant
singular vectors would exactly identify the two diagonal blocks.
Once the noise entries are inserted, however,
the dominant left singular vector
of the above matrix $\mb{A}$ accurate to
two decimal places is $[.45,.37, .37, .40, .40, .43]^T$.  In other
words, there is no separation at all between the rows numbered
$1,2,3$ and those numbered $4,5,6$, so
no submatrix is identified.  

Armed with a preliminary understanding of the convex relaxed formulation, we
now present and motivate an exact (nonconvex) formulation of LAROS,
which is as follows.
$$
\ba{rll}
\min & \norm{\bX}_* + \theta\card{{\cal I}}\card{{\cal J}}\\
\st & \inprod{\mb{A}}{\bX}\geq 1,\\
\quad & x_{ij} =0, &\forall\,(i,j)\notin{\cal I}\times{\cal J},
\ea
$$
where ${\cal I}\subset\{1,\ldots,m\}$ and ${\cal J}\subset\{1,\ldots,n\}$
are unknowns (as well as $\mb{X}$). For
fixed $\cal I$ and $\cal J$, the optimal solution would be the rank-one
approximation of the submatrix $\mb{A}({\cal I},{\cal J})$ given by the
SVD, as explained above,
and the optimal
value is $\norm{\mb{A}({\cal I},{\cal J})}^{-1}+\theta \card{{\cal
    I}}\card{{\cal J}}$. The first term of the optimal value is
$\norm{\mb{A}({\cal I},{\cal J})}^{-1}$; therefore, for appropriate
selection of $\theta$, a large submatrix (in terms of $2$-norm) will
be selected. The second term is the size of the submatrix, which is a
nonconvex function. Thus, the two terms balance the twin objectives
of selecting a submatrix
with a large first singular value and
selecting a submatrix that has a relatively small number of entries.

This now motivates \refs{eq:LAROSrelaxation}: we relax the above
nonconvex formulation by replacing the cardinality term in the
objective function with the $\ell_1$-norm.
The relaxed term $\theta\Vert
\mb{X}\Vert_1$ in the objective function has the well-known effect of
favoring sparser matrices $\mb{X}$ (those with fewer nonzero entries).
Thus, the combination of the two terms seeks a low rank matrix with
many entries equal to 0.
For example, our formulation
\req{eq:LAROSrelaxation} applied
to $\mb{A}$ above identifies the 
$\{4,5,6\}\times\{4,5,6\}$ submatrix when $\theta=0.5$.  In particular,
the solution $\mb{X}$ has zeros in all positions except
$\{4,5,6\}\times\{4,5,6\}$; in these positions it has positive entries 
ranging from $0.08$ to $0.16$.

To start our more
formal analysis, 
let us consider the following general norm minimization problem.
\be
\label{eq:gnorm}
\ba{rl}
\min & \tnorm{\bX}\\
\st & \inprod{\mb{A}}{\bX}\geq 1,
\ea
\ee
where $\tnorm{\,\cdot\,}$ is an arbitrary norm function on $\R^{m\times n}$. 
(For example, the objective function $\Vert\mb{X}\Vert_*+\theta\Vert
\mb{X}\Vert_1$ appearing in \req{eq:LAROSrelaxation} is a norm.)
Using the associated dual norm $\tnorm{\,\cdot\,}^*$, 
we can relate Problem \refs{eq:gnorm} to an equivalent problem as follows.
\begin{lemma}
\label{lem:dgnorm}
Consider $\mb{A}\neq\mb{0}$. Matrix $\bX^*$ is an optimal solution of Problem \refs{eq:gnorm} if and only if $\bY^*=\tnorm{\mb{A}}^*\bX^*$ is an optimal solution of the following problem:
\be
\label{eq:dgnorm}
\ba{rl}
\max & \inprod{\mb{A}}{\bY}\\
\st & \tnorm{\bY}\leq 1.
\ea
\ee
\end{lemma}

\begin{pf}
Let $\bX^*$ be an optimal solution of Problem \refs{eq:gnorm}. Clearly, $\bX^*\neq\mb{0}$ and $\inprod{\mb{A}}{\bX^*}=1$. Apply the norm inequality, we have
$$
\tnorm{\bX^*}\cdot \tnorm{\mb{A}}^*\geq \inprod{\mb{A}}{\bX^*}=1\,\Leftrightarrow\,\tnorm{\bX^*}\geq\frac{1}{\tnorm{\mb{A}}^*}.
$$
According to Boyd and Vandenberghe \cite{Boyd04}, the dual of the dual norm is the original norm and the norm inequality is tight: for any $\mb{A}$, there is always an $\bX\neq\mb{0}$ such that the equality holds (for finite-dimensional vector spaces). Since $\bX^*$ is an optimal solution of Problem \refs{eq:gnorm},
$$
\tnorm{\bX^*}=\frac{1}{\tnorm{\mb{A}}^*}.
$$

Let $\bY^*=\tnorm{\mb{A}}^*\bX^*$, we then have: $\tnorm{\bY^*}=1$ and $\inprod{\mb{A}}{\bY^*}=\tnorm{\mb{A}}^*$. We also have:
$$
\ba{rl}
\tnorm{\mb{A}}^*=\max & \inprod{\mb{A}}{\bY}\\
\st & \tnorm{\bY}\leq 1.
\ea
$$
Thus $\bY^*$ is indeed an optimal solution of Problem \refs{eq:dgnorm}. Using similar arguments, we can prove that conversely, if $\bY^*$ is an optimal solution of Problem \refs{eq:dgnorm}, then $\bX^*=(\tnorm{\mb{A}}^*)^{-1}\bY^*$ is an optimal solution of Problem \refs{eq:dgnorm}.
\end{pf}

The next lemma characterizes the set of all optimal solutions of Problem \refs{eq:dgnorm}.

\begin{lemma}
\label{lem:sgd}
The set of all optimal solutions of Problem \refs{eq:dgnorm} with $\mb{A}\neq\mb{0}$ is the subgradient of the dual norm function $\tnorm{\,\cdot\,}^*$ at $\mb{A}$, $\partial\tnorm{\mb{A}}^*$.
\end{lemma}

\begin{pf}
Let $\bY^*$ be an optimal solution of Problem \refs{eq:dgnorm}, we have $\tnorm{\bY^*}=1$ since $\mb{A}\neq\mb{0}$. Thus we have: $\tnorm{\mb{A}}^*=\inprod{\mb{A}}{\bY^*}$. For an arbitrary matrix $\mb{B}\in\R^{m\times n}$,
$$
\tnorm{\mb{A}+\mb{B}}^*\geq\inprod{\mb{A}+\mb{B}}{\bY^*}=\tnorm{\mb{A}}^* + \inprod{\mb{B}}{\bY^*}.
$$
Thus $\bY^*\in\partial\tnorm{\mb{A}}^*$.

Now consider $\mb{Y}\in\partial\tnorm{\mb{A}}^*$:
$$
\tnorm{\mb{B}}^*\geq\tnorm{\mb{A}}^* + \inprod{\mb{B}-\mb{A}}{\bY}\Leftrightarrow\inprod{\mb{A}}{\bY}-\tnorm{\mb{A}}^*\geq\inprod{\mb{B}}{\bY}-\tnorm{\mb{B}}^*,\quad\forall\,\mb{B}\in\R^{m\times n}.
$$
With $\mb{B}=\mb{0}$ and $\mb{B}=2\mb{A}$, we obtain the equality $\inprod{\mb{A}}{\bY}=\tnorm{\mb{A}}^*>0$. We have:
$$
\inprod{\mb{A}}{\bY}=\tnorm{\mb{A}}^*\leq\tnorm{\mb{A}}^*\tnorm{\bY} \Leftrightarrow (\tnorm{\bY}-1)\tnorm{\mb{A}}^*\geq 0\Rightarrow\tnorm{\bY}_*\geq 1.
$$

In addition, $\inprod{\mb{B}}{\bY}-\tnorm{\mb{B}}^*\leq 0$ for all $\mb{B}\in\R^{m\times n}$. The norm inequality $\inprod{\mb{B}}{\bY}\leq\tnorm{\mb{B}}^*\tnorm{\bY}$ is tight and $\bY\neq\mb{0}$ ($\tnorm{\bY}\geq 1$); therefore, there exists $\mb{B}\neq\mb{0}$ such that $\inprod{\mb{B}}{\bY}=\tnorm{\mb{B}}^*\tnorm{\bY}$. Thus we have:
$$
\tnorm{\mb{B}}^*\tnorm{\bY}-\tnorm{\mb{B}}^*\leq 0\Leftrightarrow(\tnorm{\bY}-1)\tnorm{\mb{B}}^*\leq 0\Rightarrow\tnorm{\bY}\leq 1.
$$
Thus $\tnorm{\bY}=1$ and $\inprod{\mb{A}}{\bY}=\tnorm{\mb{A}}^*$, the optimal value of Problem \refs{eq:dgnorm}, which means $\bY$ is an optimal solution of Problem \refs{eq:dgnorm}. 
\end{pf}

Lemma \ref{lem:dgnorm} and \ref{lem:sgd} show that the set of all optimal solutions of Problem \refs{eq:gnorm} is $(\tnorm{\mb{A}}^*)^{-1}\partial\tnorm{\mb{A}}^*$. The uniqueness of the optimal solution of Problem \refs{eq:gnorm} is equivalent to the differentiability of the dual norm function $\tnorm{\,\cdot\,}^*$ at $\mb{A}$. These results are summarized in the following theorem.

\begin{theorem}
\label{thm:optsol}
Consider $\mb{A}\neq\mb{0}$. The following statements are true:
\ben
\item[(i)] The set of optimal solutions of Problem \refs{eq:gnorm} is $(\tnorm{\mb{A}}^*)^{-1}\partial\tnorm{\mb{A}}^*$.
\item[(ii)] Problem \refs{eq:gnorm} has a unique optimal solution if and only if the dual norm function $\tnorm{\,\cdot\,}^*$ is differentiable at $\mb{A}$.
\een
\end{theorem}

If the norm is set to be the nuclear norm, we obtain the following minimization problem, which has been used 
\cite{Fazel:02, RechtFazelParrilo, CandesRecht, AmesVavasis}
as a relaxation of rank minimization optimization 
problems:
\be
\label{eq:nuclear}
\ba{rl}
\min & \norm{\bX}_*\\
\st & \inprod{\mb{A}}{\bX}\geq 1.
\ea
\ee

The dual norm of the nuclear norm is the spectral norm. According to Zi\c{e}tak \cite{Zietak93}, if $\mb{A}=\mb{U}\mb{\Sigma}\mb{V}^T$ is a singular value decomposition of $\mb{A}$ and $s$ is the multiplicity of the largest singular value of $\mb{A}$, the subgradient $\partial\norm{\mb{A}}$ is written as follows:
$$
\partial\norm{\mb{A}}=\left\{\mb{U}
\begin{bmatrix}
\mb{S} &\mb{0}\cr
\mb{0} & \mb{0}
\end{bmatrix}
\mb{V}^T: \mb{S}\in{\cal S}^s_+, \norm{\mb{S}}_*=1\right\}.
$$
Clearly, the largest rank-one approximation of $\mb{A}$, $\mb{u}_{1}\mb{v}_1^T$, always belongs to the subgradient $\partial\norm{\mb{A}}$. The description of the subgradient $\partial\norm{\mb{A}}$ shows that the maximum possible rank of an optimal solution of Problem \refs{eq:nuclear} is the multiplicity of the largest singular value of $\mb{A}$. In addition, the spectral norm function $\norm{\,\cdot\,}$ is not differentiable in general. The uniqueness of the optimal solution of Problem \refs{eq:nuclear} is equivalent to the differentiability of the spectral norm function $\norm{\,\cdot\,}$ at $\mb{A}$. The necessary and sufficient condition is $s=1$ or equivalently, $\sigma_1(\mb{A})>\sigma_2(\mb{A})$. In the case of unique optimal solution, we obtain the largest rank-one approximation of $\mb{A}$ (up to the scaling factor $\norm{\mb{A}}^{-1}$). These results are stated in the following proposition:
\begin{proposition}
\label{prop:nuclear}
Consider $\mb{A}\neq\mb{0}$. The following statements are true:
\ben
\item[(i)] The set of optimal solutions of Problem \refs{eq:nuclear} is $\norm{\mb{A}}^{-1}\partial\norm{\mb{A}}$.
\item[(ii)] The largest rank-one approximation of $\mb{A}$ is an optimal solution of Problem \refs{eq:nuclear} and it is the unique solution if and only if $\sigma_1(\mb{A})>\sigma_2(\mb{A})$.
\een
\end{proposition}  

Similar to low-rank minimization problems with nuclear norm approximation, sparse optimization problems can be approximately handled by the (vector) $\ell_1$-norm function $\norm{\,\cdot\,}_1$. Let us consider the following problem
\be
\label{eq:l1norm}
\ba{rl}
\min & \norm{\bX}_1\\
\st & \inprod{\mb{A}}{\bX}\geq 1.
\ea
\ee
The dual norm of $\ell_1$-norm is the (vector) infinity norm 
$\norm{\,\cdot\,}_{\infty}$, i.e., the maximum absolute entry of the matrix,
and the subgradient $\partial\norm{\mb{A}}_{\infty}$ can be written as follows,
$$
\partial\norm{\mb{A}}_{\infty}=\mbox{conv}\left\{\mbox{sgn}(a_{ij})\mb{E}_{ij}\mid (i,j)\in\arg\max_{(k,l)}\card{a_{kl}}\right\},
$$
where $\mb{E}_{ij}$ is the unit matrix in $\R^{m\times n}$ with $\mb{E}_{ij}(i,j)=1$. The sparsity of the optimal solution $\bX^*$ of Problem \refs{eq:l1norm} is clearly related to the multiplicity of the maximum absolute value of elements of $\mb{A}$. Applying Theorem \ref{thm:optsol} for this particular $\ell_1$-norm, we obtain the following results:
\begin{proposition}
\label{prop:l1norm}
Consider $\mb{A}\neq\mb{0}$. The following statements are true:
\ben
\item[(i)] The set of optimal solutions of Problem \refs{eq:l1norm} is $\norm{\mb{A}}_{\infty}^{-1}\partial\norm{\mb{A}}_{\infty}$.
\item[(ii)] The matrix $\mathop{\rm sgn}(a_{ij})\mb{E}_{ij}$, where $\displaystyle (i,j)\in\arg\max_{(k,l)}\card{a_{kl}}$, and $\mathop{\rm sgn}(\cdot)$ is the usual sign function, is an optimal solution of Problem \refs{eq:l1norm} and it is the unique solution if and only if $\card{a_{ij}}>\card{a_{kl}}$ for all $(k,l)\neq (i,j)$.
\een
\end{proposition} 

As mentioned above, finding a low-rank submatrix clearly involves both low-rank and sparse optimization (with a specific sparse structure). Let us 
return to the parametric optimization problem 
\req{eq:LAROSrelaxation}
proposed at the beginning of this section
$$
\ba{rl}
\min & \norm{\bX}_*+\theta\norm{\bX}_1\\
\st & \inprod{\mb{A}}{\bX}\geq 1,
\ea
$$
where $\theta\geq 0$. Clearly, if $\theta=0$, we obtain Problem \refs{eq:nuclear} and when $\theta\rightarrow\infty$, we approach Problem \refs{eq:l1norm}. This optimization problem clearly addresses both low-rank and sparse requirements of the solution $\bX$. We now would like to characterize the set of optimal solutions of the problem.

The objective function $\norm{\bX}_*+\theta\norm{\bX}_1$ is a norm function since $\theta\geq 0$. Denote $\norm{\bX}_{\theta}$ to be this parametric norm of $\bX$, 
$$
\norm{\bX}_{\theta}:=\norm{\bX}_*+\theta\norm{\bX}_1
$$
and consider its dual norm function $\norm{\,\cdot\,}_{\theta}^*$.  Clearly, Problem \refs{eq:LAROSrelaxation} is a special case of Problem \refs{eq:gnorm}. The set of optimal solutions of Problem \req{eq:LAROSrelaxation} can therefore be characterized as follows:
\begin{proposition}
\label{thm:msol}
Consider $\mb{A}\neq\mb{0}$. The following statements are true:
\ben
\item[(i)] The set of optimal solutions of Problem \refs{eq:LAROSrelaxation} is $(\norm{\mb{A}}_{\theta}^*)^{-1}\partial\norm{\mb{A}}_{\theta}^*$.
\item[(ii)] There is a unique solution if and only if the dual norm function $\norm{\,\cdot\,}_{\theta}^*$ is differentiable at $\mb{A}$.
\een
\end{proposition}

We now focus on deriving some properties of the dual norm $\norm{\,\cdot\,}_{\theta}^*$. We have:
\be
\label{eq:dnorm}
\ba{rl}
\norm{\mb{A}}_{\theta}^*=\max & \inprod{\mb{A}}{\bX}\\
\st & \norm{\bX}_{\theta}\leq 1.
\ea
\ee
We will use the gauge function and its dual polar function (see Rockafellar \cite{Rockafellar70} for more details) to compute this dual norm.
\begin{proposition}
\label{prop:dnorm}
The dual norm $\norm{\mb{A}}_{\theta}^*$ with $\theta>0$ is the optimal value of the following optimization problem:
\be
\label{eq:gfunc}
\ba{rl}
\norm{\mb{A}}_{\theta}^*=\min & \max\left\{\norm{\bY},\theta^{-1}\norm{\bZ}_{\infty}\right\}\\
\st & \bY + \bZ = \mb{A}.
\ea
\ee
\end{proposition}

\begin{pf}
Consider the closed unit ball ${\cal C}_*=\{\bX\in\R^{m\times n}\mid\norm{\bX}_*\leq 1\}$ with respect to the nuclear norm and similarly, the unit ball ${\cal C}_1$ with respect to the $\ell_1$-norm $\norm{\,\cdot\,}_1$. We have the polar of ${\cal C}_*$ is the closed unit ball with respect to the spectral norm, ${\cal C}_*^\circ={\cal C}$. Similarly, we have: ${\cal C}_1^\circ = {\cal C}_{\infty}$, the unit ball with respect to the infinity norm.

Using the definition of gauge functions, we have: $\norm{\bX}_* = \gamma_{{\cal C}_*}(\bX)=\min\{\lambda\geq 0\mid\bX\in\lambda{\cal C}_*\}$. In addition, the support function $\sigma_{{\cal S}}(\bX)=\max\{\inprod{\bX}{\bY}\mid\bY\in{\cal S}\}$ is the gauge function of ${\cal S}^\circ$ for all symmetric closed bounded convex set with $\mb{0}\in\mbox{int}({\cal S})$. All unit balls satisfy these conditions; therefore, we obtain the well-known results $\norm{\bX}_*=\gamma_{{\cal C}_*}(\bX)=\sigma_{{\cal C}_*^\circ}(\bX)=\sigma_{{\cal C}}(\bX)$ and $\norm{\bX}_1=\gamma_{{\cal C}_1}(\bX)=\sigma_{{\cal C}_1^\circ}(\bX)=\sigma_{{\cal C}_{\infty}}(\bX)$.

Now consider the unit ball ${\cal C}_{\theta}=\{\bX\in\R^{m\times n}\mid\norm{\bX}_*+\theta\norm{\bX}_1\leq 1\}$, we have:
$$
{\cal C}_{\theta}=\{\bX\in\R^{m\times n}\mid\sigma_{{\cal C}}(\bX)+\theta\sigma_{{\cal C}_{\infty}}(\bX)\leq 1\}.
$$

Applying the definition of support functions, we have: $\sigma_{{\cal C}}(\bX)+\theta\sigma_{{\cal C}_{\infty}}(\bX)=\sigma_{{\cal C}+\theta{\cal C}_{\infty}}(\bX)$, where ${\cal C}+\theta{\cal C}_{\infty}$ is the Minkowski sum of two sets, ${\cal C}$ and $\theta{\cal C}_{\infty}$. This set satisfies all the conditions above; therefore, $\sigma_{{\cal C}+\theta{\cal C}_{\infty}}(\bX)=\gamma_{({\cal C}+\theta{\cal C}_{\infty})^\circ}(\bX)$. Thus
$$
{\cal C}_{\theta}=\{\bX\in\R^{m\times n}\mid \gamma_{({\cal C}+\theta{\cal C}_{\infty})^\circ}(\bX)\leq 1\}=({\cal C}+\theta{\cal C}_{\infty})^\circ.
$$

We also have: $\norm{\mb{A}}_{\theta}^*=\sigma_{{\cal C}_{\theta}}(\mb{A})$. Thus
$$
\norm{\mb{A}}_{\theta}^*=\gamma_{{\cal C}_{\theta}^\circ}(\mb{A})=\gamma_{{\cal C}+\theta{\cal C}_{\infty}}(\mb{A}).
$$

We have: $\gamma_{{\cal C}+\theta{\cal C}_{\infty}}(\mb{A})=\min\{\lambda\geq 0\mid \mb{A}\in \lambda({\cal C}+\theta{\cal C}_{\infty})\}$ or equivalently,
$$
\ba{rl}
\gamma_{{\cal C}+\theta{\cal C}_{\infty}}(\mb{A})=\min & \lambda\\
\st & \mb{A}=\mb{Y}+\mb{Z},\\
\quad & \norm{\bY}\leq \lambda,\\
\quad &\theta^{-1}\norm{\bZ}_{\infty}\leq \lambda,\\
\quad &\lambda\geq 0.
\ea
$$

Rewriting the minimization problem above, we obtain the final result as shown in \refs{eq:gfunc}:
$$
\ba{rl}
\norm{\mb{A}}_{\theta}^*=\min & \max\left\{\norm{\bY},\theta^{-1}\norm{\bZ}_{\infty}\right\}\\
\st & \bY + \bZ = \mb{A}.
\ea
$$
\end{pf}

We can now derive the optimality conditions for both problems \refs{eq:dnorm} and \refs{eq:gfunc}:
\begin{lemma}
\label{lem:optimality}
Nonzero feasible solutions $\bX$ and $(\bY,\bZ)$ are optimal for Problem \refs{eq:dnorm} and \refs{eq:gfunc} respectively if and only if they satisfy the conditions below:
\ben
\item[(i)] $\norm{\bY}=\theta^{-1}\norm{\bZ}_{\infty}$,
\item[(ii)] $\bX\in\alpha\partial\norm{\bY}$, $\alpha\geq 0$,
\item[(iii)] $\bX\in\beta\partial\norm{\bZ}_{\infty}$, $\beta\geq 0$, and
\item[(iv)] $\alpha+\theta\beta=1$.
\een
\end{lemma}  

\begin{pf}
We first prove the weak duality result. Consider feasible solutions $\bX$ and $(\bY,\bZ)$ for Problem \refs{eq:dnorm} and \refs{eq:gfunc} respectively, we have:
$$
\ba{rcl}
\inprod{\bX}{\mb{A}} &=& \inprod{\bX}{\bY}+\inprod{\bX}{\bZ}\\
\quad &\leq& \norm{\bX}_*\norm{\bY} + \norm{\bX}_1\norm{\bZ}_{\infty}\\
\quad &\leq& \norm{\bX}_*\max\{\norm{\bY},\theta^{-1}\norm{\bZ}_{\infty}\} + \theta\norm{\bX}_1\max\{\norm{\bY},\theta^{-1}\norm{\bZ}_{\infty}\}\\
\quad &=&(\norm{\bX}_*+\theta\norm{\bX}_1)\max\{\norm{\bY},\theta^{-1}\norm{\bZ}_{\infty}\}\\
\quad &\leq& \max\{\norm{\bY},\theta^{-1}\norm{\bZ}_{\infty}\}. 
\ea
$$

The strong duality result shows that $\bX$ and $(\bY,\bZ)$ are the optimal solutions if and only if $\inprod{\bX}{\mb{A}}=\max\{\norm{\bY},\theta^{-1}\norm{\bZ}_{\infty}\}$. This happens if and only if all the conditions below are satisfied:
\ben
\item[(i)] $\inprod{\bX}{\bY}=\norm{\bX}_*\norm{\bY}$ and $\inprod{\bX}{\bZ}=\norm{\bX}_1\norm{\bZ}_{\infty}$,
\item[(ii)] $\norm{\bY}=\max\{\norm{\bY},\theta^{-1}\norm{\bZ}_{\infty}\}=\theta^{-1}\norm{\bZ}_{\infty}$, and
\item[(iii)] $\norm{\bX}_*+\theta\norm{\bX}_1=1$.
\een

The first two conditions are equivalent to the fact that $\bX=\alpha\partial\norm{\bY}$, where $\alpha=\norm{\bX}_*$, and $\bX=\beta\partial\norm{\bZ}_{\infty}$, where $\beta=\norm{\bX}_1$. The second condition is simply $\norm{\bY}=\theta^{-1}\norm{\bZ}_{\infty}$ and the third condition is equivalent to $\alpha + \theta\beta=1$. Thus we have proved the necessary and sufficient optimality conditions for Problem \refs{eq:dnorm} and \refs{eq:gfunc}.
\end{pf}

Using these optimality conditions, we can obtain simple sufficient conditions for the uniqueness of the optimal solution $\bX$:
\begin{proposition}
\label{prop:unique}
Consider the feasible solution $\bX$ of Problem \refs{eq:dnorm}. If there exists $(\bY,\bZ)$ that satisfies the conditions below,
\ben
\item[(i)] $\bY+\bZ=\mb{A}$ and $\norm{\bY}=\theta^{-1}\norm{\bZ}_{\infty}$,
\item[(ii)] $\bX\in\alpha\partial\norm{\bY}$, $\alpha\geq 0$,
\item[(iii)] $\bX\in\beta\partial\norm{\bZ}_{\infty}$, $\beta\geq 0$,
\item[(iv)] $\alpha+\theta\beta=1$, and
\item[(v)] $\norm{\,\cdot\,}$ is differentiable at $\bY$ or $\norm{\,\cdot\,}_{\infty}$ is differentiable at $\bZ$,
\een
then $\bX$ is the unique optimal solution of Problem \refs{eq:dnorm}.
\end{proposition}

\begin{pf}
Using the first four conditions, we can prove that $\bX$ is an optimal solution of Problem \refs{eq:dnorm} and $(\bY,\bZ)$ is an optimal solution of Problem \refs{eq:gfunc}. Now assume that $\norm{\,\cdot\,}$ is differentiable at $\bY$, we have: $\partial\norm{\bY}$ is a singleton, $\partial\norm{\bY}=\{\mb{V}\}$. Thus we have:
$$
\norm{\bX}_1=\alpha\norm{\mb{V}}_1=\beta\Rightarrow \alpha(1+\theta\norm{\mb{V}}_1)=0\Rightarrow \alpha=\frac{1}{1+\theta\norm{\mb{V}}_1}.
$$
Assume there is another optimal solution $\bar{\bX}\neq\bX$ of Problem \refs{eq:dnorm}. Applying Lemma \ref{lem:optimality}, we will have $\bar{\bX}\in\bar{\alpha}\partial\norm{\bY}$ and similarly $\bar{\bX}\in\bar{\beta}\partial\norm{\bZ}_{\infty}$ with $\bar{\alpha}+\theta\bar{\beta}=1$. Same calculation results in $\bar{\alpha}=\alpha$ (contradiction). Thus $\bX$ is the unique optimal solution of Problem \refs{eq:dnorm}. Similar arguments can be used to prove the uniqueness of $\bX$ if $\norm{\,\cdot\,}_{\infty}$ is differentiable at $\bZ$.
\end{pf}

Proposition \ref{prop:unique} relies on dual solutions $\bY$ and $\bZ$ to show the uniqueness of the primal solution $\bX$. Next, we will focus on the low-rank and sparse property of the optimal solution $\bX$ for different values of $\theta$. The following theorem provides the sufficient conditions on matrix $\mb{A}$ for the rank-one property (and uniqueness) of the optimal solution $\bX$ when $\theta$ is small enough.
\begin{theorem}
\label{theorem:unique}
If $\mb{A}$ satisfies the condition $\sigma_1(\mb{A})>\sigma_2(\mb{A})$, then Problem \refs{eq:dnorm} has a (unique) rank-one optimal solution $\bX$ for all $0\leq\theta<\theta_A$, where $\displaystyle\theta_A=\frac{1}{\sqrt{mn}}\left(\frac{\sigma_1(\mb{A})-\sigma_2(\mb{A})}{3\sigma_1(\mb{A})-\sigma_2(\mb{A})}\right)$.
\end{theorem}

\begin{pf}
The optimality conditions in Lemma \ref{lem:optimality} show that there exist $\bY$ and $\bZ$ such that $\mb{A}=\bY + \bZ$, $\norm{\bZ}_{\infty}=\theta\norm{\bY}$ and $\bX\in\alpha\partial\norm{\bY}$. Applying 
a standard perturbation theorem of singular values
(see Cor.~8.6.2 of \cite{GVL}), we have:
$$
\card{\sigma_i(\mb{A})-\sigma_i(\bY)}\leq \norm{\bZ},\quad i=1,2.
$$
We also have: $\norm{\bZ}\leq\sqrt{mn}\norm{\bZ}_{\infty}$. Thus
$$
\norm{\bZ}\leq\sqrt{mn}\left(\theta\norm{\bY}\right)=\sqrt{mn}\left(\theta\sigma_1(\bY)\right).
$$

For all $0\leq\theta<\theta_A$, we have:
$$
\sigma_1(\bY)\leq\sigma_1(\mb{A})+\norm{\bZ}\leq\sigma_1(\mb{A})+\sqrt{mn}\left(\theta\sigma_1(\bY)\right)<\sigma_1(\mb{A})+\sqrt{mn}\left(\theta_A\sigma_1(\bY)\right).
$$
This implies
$$
(1-\theta_A\sqrt{mn})\sigma_1(\bY)<\sigma_1(\mb{A})\Leftrightarrow\frac{2\sigma_1(\mb{A})}{3\sigma_1(\mb{A})-\sigma_2(\mb{A})}\sigma_1(\bY)<\sigma_1(\mb{A})\Leftrightarrow \sigma_1(\bY)<\frac{1}{2}(3\sigma_1(\mb{A})-\sigma_2(\mb{A})).
$$
We then have:
$$
\norm{\bZ}\leq \sqrt{mn}\left(\theta\sigma_1(\bY)\right)<\sqrt{mn}\left(\theta_A\sigma_1(\bY)\right)<\frac{1}{2}(\sigma_1(\mb{A})-\sigma_2(\mb{A})).
$$
Thus
$$
\sigma_1(\bY)\geq\sigma_1(\mb{A})-\norm{\bZ}>\frac{1}{2}(\sigma_1(\mb{A})+\sigma_2(\mb{A}))>\sigma_2(\mb{A})+\norm{\bZ}\geq\sigma_2(\bY).
$$

We have $\sigma_1(\bY)>\sigma_2(\bY)$; therefore, $\norm{\,\cdot\,}$ is differentiable at $\bY$. According to Proposition \ref{prop:unique}, we have $\bX$ is the unique rank-one optimal solution of Problem \refs{eq:dnorm}.
\end{pf}

The last result of this section concerns the nonnegativity of 
$\mb{X}$.  If $\mb{A}$ is nonnegative, then one might expect 
$\mb{X}$ to be nonnegative.  For $\theta=0$ or $\theta=\infty$, this
is certainly true by preceding results in this section.  It is not always
necessarily true for intermediate values of $\theta$.  The following
theorem shows that, at least in the rank-one case, nonnegativity
is assured.

\begin{theorem}
\label{theorem:nneg}
Consider the set of optimal solutions of Problem \refs{eq:dnorm}
when $\mb{A}\ge \mb{0}$.  We have:
\ben
\item[(i)] If Problem \refs{eq:dnorm} has a rank-one optimal solution, then there exists a nonnegative rank-one optimal solution.
\item[(ii)] If $\theta>1$, then all optimal solutions of Problem \refs{eq:dnorm} are nonnegative.
\een
\end{theorem}

\begin{pf}
\ben
\item[(i)] Consider a rank-one optimal solution $\bX$, $\bX=\sigma\mb{u}\mb{v}^T$, of Problem \refs{eq:dnorm}. We prove that $\card{\bX}=\sigma\card{\mb{u}}\card{\mb{v}}^T\geq\mb{0}$ is also an optimal solution. 
Let $\tilde{\bX}$ denote $\card{\bX}$.  
We have:
$$
\norm{\tilde{\bX}}_{\theta}=\norm{\tilde{\bX}}_*+\theta\norm{\tilde{\bX}}_1=\norm{\bX}_*+\theta\norm{\bX}_1.
$$
In addition, $\inprod{\mb{A}}{\tilde{\bX}}\geq\inprod{\mb{A}}{\bX}$ since $\mb{A}\geq\mb{0}$. Thus clearly $\tilde{\bX}$ is also an optimal solution.
\item[(ii)] Assume that there exists an optimal solution $\bX$ of Problem \refs{eq:dnorm} is not nonnegative. Without loss of generality, assume $x_{11}<0$. Consider $\bX(\eps)=\bX+\eps\mb{E}_{11}$, where $\eps>0$ and $\mb{E}_{11}$ is the matrix of all zeros except the element $\mb{E}_{11}(1,1)=1$, we have:
$$
\norm{\bX(\eps)}_*\leq\norm{\bX}_*+\eps\norm{\mb{E}_{11}}_*=\norm{\bX}_*+\eps.
$$
In addition, $\norm{\bX(\eps)}_1=\norm{\bX}_1-\eps$ if $\eps\leq\card{x_{11}}$. Therefore, we have:
$$
\norm{\bX(\eps)}_{\theta}\leq\norm{\bX}+(1-\theta)\eps=1+(1-\theta)\eps< 1,\quad\forall\, 0<\eps\leq\card{x_{11}}.
$$
Here we assume that $\mb{A}\neq\mb{0}$; therefore, $\norm{\bX}_{\theta}=1$. We also have
$$
\inprod{\mb{A}}{\bX(\eps)}=\inprod{\mb{A}}{\bX}+\eps a_{11}\geq\inprod{\mb{A}}{\bX}.
$$
Now consider $\displaystyle\bar{\bX}=\frac{1}{1+(1-\theta)\eps}\bX(\eps)$. Clearly, $\norm{\bar{\bX}}_{\theta}=1$ and $\inprod{\mb{A}}{\bar{\bX}}>\inprod{\mb{A}}{\bX}>0$ (contradiction). Thus all optimal solutions of Problem \refs{eq:dnorm} are nonnegative if $\theta>1$.
\een
\end{pf}

\section{Sparsity}
\label{sec:sparsity}

As mentioned in the introduction, the penalty term
$\theta\Vert\mb{X}\Vert_1$ in the objective function
of \req{eq:LAROSrelaxation} is intended to promote sparsity
of $\mb{X}$.  For some very simple convex optimization problems
with an $\ell_1$ penalty term, e.g., the unconstrained problem
of minimizing $\Vert \vect{x}-\vect{c}\Vert_2+\theta\Vert\vect{x}\Vert_1$ 
for a given
vector $\vect{c}$, it is
known that sparsity increases monotonically with $\theta$
(i.e., if $\vect{x}^*_1$ is the optimizer for $\theta_1$ and
$\vect{x}^*_2$ is the optimizer for $\theta_2$ with $\theta_1\le \theta_2$,
then the indices of nonzeros of $\vect{x}_2^*$ are a subset of the indices
of nonzeros of $\vect{x}_1^*$).

For a more complicated problem such as \req{eq:LAROSrelaxation}, monotonicity
does not hold in general.  But nonetheless, some weaker statements
about the relationship between $\theta$ and sparsity are possible.
Two such results are derived in this section.
We start with a lemma that leads to a sparsity result.

\begin{lemma}
\label{lemma:roptimality}
Assume $\bX=\sigma\mb{u}\mb{v}^T$, where $\norm{\mb{u}}_2=\norm{\mb{v}}_2=1$, $\mb{u}\geq\mb{0}$, and $\mb{v}\geq\mb{0}$, is the optimal solution of Problem \refs{eq:dnorm}. If $u_i>u_j=0$ then
\be
\label{eq:ropt}
\norm{A}_{\theta}^*=\frac{\mb{a}_i^T\mb{v}}{\theta\norm{\mb{v}}_1+u_i}\geq\frac{\mb{a}_j^T\mb{v}}{\theta\norm{\mb{v}}_1},
\ee
where $\mb{a}_i$ and $\mb{a}_j$ are the $i$th and $j$th row of $\mb{A}$.
\end{lemma}

\begin{pf}
We again assume here $\mb{A}\neq\mb{0}$, which means $\norm{\mb{X}}_{\theta}=1$. We have: 
$$
\norm{\bX}_{\theta}=\sigma + \theta\sigma\norm{\mb{u}}_1\norm{\mb{v}}_1=1\Leftrightarrow\sigma=\frac{1}{1+\theta\norm{\mb{u}}_1\norm{\mb{v}}_1}.
$$ 
Consider $\bX(\eps)=\sigma(\mb{u}+\eps\mb{e}_i)\mb{v}^T$, where $\eps\geq 0$ and $\mb{e}_i$ is the $i$th unit vector, we have:
$$
\norm{\mb{u}+\eps\mb{e}_i}_2^2=\norm{\mb{u}}_2^2+(u_i+\eps)^2-u_i^2=1+2u_i\eps+\eps^2.
$$
Thus $\norm{\bX(\eps)}_*=\sigma\sqrt{1+2u_i\eps+\eps^2}$. On the other hand, $\norm{\bX(\eps)}_1=\sigma(\norm{\mb{u}}_1+\eps)\norm{\mb{v}}_1$. So we have:
$$
\norm{\bX(\eps)}_{\theta}=1+\sigma\left(\sqrt{1+2u_i\eps+\eps^2}-1+\theta\eps\norm{\mb{v}}_1\right)=1+\sigma\eps\left(\frac{2u_i+\eps}{\sqrt{1+2u_i\eps+\eps^2}+1}+\theta\norm{\mb{v}}_1\right).
$$

Let $\displaystyle\alpha = \frac{2u_i+\eps}{\sqrt{1+2u_i\eps+\eps^2}+1}+\theta\norm{\mb{v}}_1$ and consider $\displaystyle\bar{\bX}=\frac{\bX(\eps)}{1+\sigma\eps\alpha}$, we have: $\norm{\bar{\bX}}_{\theta}=1$ and
$$
\inprod{\mb{A}}{\bar{\bX}}=\frac{\inprod{\mb{A}}{\bX(\eps)}}{1+\sigma\eps\alpha}=\frac{\norm{\mb{A}}_{\theta}^*+\sigma\eps\mb{a}_i^T\mb{v}}{1+\sigma\eps\alpha}\leq \norm{\mb{A}}_{\theta}^*.
$$
With $\sigma>0$ and $\eps>0$, we obtain the following inequality
$$
\mb{a}_i^T\mb{v}\leq\left(\frac{2u_i+\eps}{\sqrt{1+2u_i\eps+\eps^2}+1}+\theta\norm{\mb{v}}_1\right)\norm{\mb{A}}_{\theta}^*.
$$
Taking the limit $\eps\rightarrow 0^+$, we have:
$$
\norm{\mb{A}}_{\theta}^*\geq\frac{\mb{a}_i^T\mb{v}}{\theta\norm{\mb{v}}_1+u_i}.
$$

Now consider the case in which $u_i>0$ and set $\bX(\eps)=\sigma(\mb{u}-\eps\mb{e}_i)\mb{v}^T$, where $0\leq\eps\leq u_i$. Similarly, we have:
$$
\norm{\bX(\eps)}_{\theta}=1-\sigma\eps\left(\frac{2u_i-\eps}{\sqrt{1-2u_i\eps+\eps^2}+1}+\theta\norm{\mb{v}}_1\right).
$$
This implies the following inequality
$$
\mb{a}_i^T\mb{v}\geq\left(\frac{2u_i-\eps}{\sqrt{1-2u_i\eps+\eps^2}+1}+\theta\norm{\mb{v}}_1\right)\norm{\mb{A}}_{\theta}^*.
$$
Again, taking the limit $\eps\rightarrow 0^+$, we have:
$$
\norm{\mb{A}}_{\theta}^*\leq\frac{\mb{a}_i^T\mb{v}}{\theta\norm{\mb{v}}_1+u_i}.
$$
{From} these two results, we can see that if $u_i>u_j=0$, then
$$
\norm{A}_{\theta}^*=\frac{\mb{a}_i^T\mb{v}}{\theta\norm{\mb{v}}_1+u_i}\geq\frac{\mb{a}_j^T\mb{v}}{\theta\norm{\mb{v}}_1}.
$$
\end{pf}

Since the roles of columns and rows are interchangeable, we also have the following result. If $v_k>v_l=0$, then
\be
\label{eq:copt}
\norm{A}_{\theta}^*=\frac{\mb{u}^T\mb{A}_k}{\theta\norm{\mb{u}}_1+v_k}\geq\frac{\mb{u}^T\mb{A}_l}{\theta\norm{\mb{u}}_1},
\ee
where $\mb{A}_k$ and $\mb{A}_l$ are $k$-th and $l$-th column of $\mb{A}$.

The sparsity structure of $\bX=\sigma\mb{u}\mb{v}^T$ depends on the sparsity structure of $\mb{u}$ and $\mb{v}$. The results obtained above help us derive some conditions under which a row (or column) of $\bX$ is zero.

\begin{corollary}
\label{col:zrow}
Consider two rows $\mb{a}_i^T$ and $\mb{a}_j^T$ of matrix $\mb{A}$. If $\min_{k}\left\{a_{ik}\right\}\geq\alpha\max_k\left\{a_{jk}\right\}$, where $\alpha>1$, then for every $\displaystyle\theta>\frac{1}{\alpha-1}$ and for every nonnegative rank-one optimal solution $\bX$ of Problem \refs{eq:dnorm}, the $j$th row of $\bX$ is zero.
\end{corollary}

\begin{pf}
Assume that $\bX=\sigma\mb{u}\mb{v}^T$, where $\mb{u}\geq\mb{0}$ and $\mb{v}\geq\mb{0}$, and $\norm{\mb{u}}_2=\norm{\mb{v}}_2=1$. We have:
$$
\frac{\mb{a}_i^T\mb{v}}{\theta\norm{\mb{v}}_1+u_i}\geq\frac{\min_k\{a_{ik}\}\norm{\mb{v}}_1}{\theta\norm{\mb{v}}_1+1}\geq\frac{\alpha\max_k\{a_{jk}\}\norm{\mb{v}}_1}{\theta\norm{\mb{v}}_1+1}=\frac{\alpha\max_k\{a_{jk}\}}{\displaystyle\theta+\frac{1}{\norm{\mb{v}}_1}}\geq\frac{\alpha\max_k\{a_{jk}\}}{\theta+1},
$$
since $\norm{\mb{v}}_1\geq\norm{\mb{v}}_2=1$. On the other hand, we also have the following inequality
$$
\frac{\mb{a}_j^T\mb{v}}{\theta\norm{\mb{v}}_1+u_j}\leq\frac{\max_k\{a_{jk}\}\norm{\mb{v}}_1}{\theta\norm{\mb{v}}_1}=\frac{\max_k\{a_{jk}\}}{\theta}.
$$
We have:
$$
\frac{\alpha\max_k\{a_{jk}\}}{\theta+1}-\frac{\max_k\{a_{jk}\}}{\theta}=\max_k\{a_{kj}\}\frac{(\alpha-1)\theta-1}{\theta(\theta+1)}>0,\quad\forall\,\theta>\frac{1}{\alpha-1}.
$$
Thus we have:
$$
\frac{\mb{a}_i^T\mb{v}}{\theta\norm{\mb{v}}_1+u_i}>\frac{\mb{a}_j^T\mb{v}}{\theta\norm{\mb{v}}_1+u_j},\quad\forall\,\theta>\frac{1}{\alpha-1},
$$
which means $u_j=0$ according to Lemma \ref{lemma:roptimality}. Thus the $j$-th row of $\bX$ is zero.
\end{pf}

We would like to use these results to build up results for columns and rows simultaneously. More exactly, consider a subset ${\cal I}\subset\{1,\ldots,m\}$ and ${\cal J}\subset\{1,\ldots,n\}$, we would like to obtain conditions on magnitudes of elements of $\mb{A}({\cal I},{\cal J})$ as compared to those of the remaining elements of $\mb{A}$ to guarantee that all rows and columns that are not in $\cal I$ and $\cal J$ have to be zero in the nonnegative rank-one optimal matrix $\bX$ of Problem \refs{eq:dnorm} for $\theta\geq\theta_0$. One of the difficulties here is that under these conditions, there is a coupling relationship between rows and columns. More exactly, in order to prove the rows that are not in $\cal I$ are zero, we need to prove the columns that are not in $\cal J$ are small or zero at the same time.

Lemma \ref{lemma:roptimality} and Corollary \ref{col:zrow} are based on local optimality conditions with respect to rows or columns. We can obtain additional results on the sparsity of the optimal solution $\bX$ using the global optimality conditions.

The following theorem states that if the weight of nonnegative matrix
$\mb{A}$ is concentrated in a
particular subblock then for $\theta$ sufficiently large, the optimal solution
$\mb{X}$ will have nonzero entries only in that subblock.  ``Concentration
of weight'' in this sense means that the average of those entries dominates
all the other entries of the matrix.  

As a special case, this theorem implies that if the maximum entry of
$\mb{A}$ is unique, then for $\theta$ sufficiently large, $\mb{X}$ will
be a singleton matrix whose unique nonzero entry corresponds to the
maximum entry of $\mb{A}$.

\begin{theorem}
\label{thm:zblock}
Assume $\mb{A}\ge\mb{0}$.
Let ${\cal I}$ and ${\cal J}$ be subsets of $\{1,\ldots,m\}$ and $\{1,\ldots,n\}$, respectively; $\card{{\cal I}}=M$ and $\card{{\cal J}}=N$. Define $\displaystyle\bar{a}({\cal I},{\cal J})=\frac{1}{MN}\sum_{(i,j)\in({\cal I},{\cal J})}a_{ij}$ and $\displaystyle a_{\max}(\overline{{\cal I},{\cal J}})=\max_{(i,j)\not\in({\cal I},{\cal J})}a_{ij}$. If $\bar{a}({\cal I},{\cal J})>a_{\max}(\overline{{\cal I},{\cal J}})$ then all optimal solutions $\bX$ of Problem \refs{eq:dnorm} are sparse, $x_{ij}=0$ for all $(i,j)\not\in({\cal I},{\cal J})$, for all $\theta>\theta_B$, where $\displaystyle\theta_B=\frac{1}{\sqrt{MN}}\left(\frac{\bar{a}({\cal I},{\cal J})\sqrt{MM}+a_{\max}(\overline{{\cal I},{\cal J}})}{\bar{a}({\cal I},{\cal J})-a_{\max}(\overline{{\cal I},{\cal J}})}\right)$.
\end{theorem}

\begin{pf}
Assume there exists an optimal solution $\bX$ of Problem \refs{eq:dnorm} such that $x_{ij}\neq 0$ for some $(i,j)\not\in({\cal I},{\cal J})$ when $\theta>\theta_B$. We have: $\theta_B>1$; therefore, according to Theorem \ref{theorem:nneg}, $\bX\geq\mb{0}$. Thus $x_{ij}>0$. We also have: $\mb{A}\neq \mb{0}$; therefore $\bX\neq\mb{0}$ and $\norm{\bX}_{\theta}=1$. Consider two cases, $a_{ij}=0$ and $a_{ij}>0$.

If $a_{ij}=0$, let $\bX_0=\bX-x_{ij}\mb{E}_{ij}$, where $\mb{E}_{ij}$ is the matrix of all zeros but the element $\mb{E}_{ij}(i,j)=1$. We have:
$\norm{\bX_0}\leq\norm{\bX}+x_{ij}$ and $\norm{\bX_0}_1=\norm{\bX}_1-x_{ij}$. Thus
$$
\norm{\bX_0}_{\theta}\leq\norm{\bX}_{\theta}+(1-\theta)x_{ij}<1,\quad\forall\,\theta>\theta_B>1.
$$
We also have: since $a_{ij}=0$, $\inprod{\mb{A}}{\bX_0}=\inprod{\mb{A}}{\bX}=\norm{\mb{A}}_{\theta}^*>0$. Thus $\bX_0\neq\mb{0}$ or $\norm{\bX_0}_{\theta}>0$. Define $\displaystyle\bX_0^s=\frac{1}{\norm{\bX_0}_{\theta}}\bX_0$, we have: $\bX_0^s$ is a feasible solution of Problem \refs{eq:dnorm} with the objective value $\displaystyle\inprod{\mb{A}}{\bX_0^s}=\frac{\norm{\mb{A}}_{\theta}^*}{\norm{\bX_0}_{\theta}}>\norm{\mb{A}}_{\theta}^*$ (contradiction).

Now consider the case when $a_{ij}>0$. Define $\mb{D}=\mb{e}_{\cal I}\mb{e}_{\cal J}^T-MNr\mb{e}_i\mb{e}_j^T$, where $\displaystyle r=\frac{\bar{a}({\cal I},{\cal J})}{a_{ij}}$, $\displaystyle\mb{e}_{\cal I}=\sum_{i\in{\cal I}}\mb{e}_i$, $\mb{e}_i\in\R^m$ is the $i$th unit vector in $\R^m$, and similarly, $\displaystyle\mb{e}_{\cal J}=\sum_{j\in{\cal J}}\mb{e}_j$, $\mb{e}_j\in\R^n$ is the $j$-th unit vector in $\R^n$. We have $r$ is well-defined since $\mb{A}>\mb{0}$ and $r>1$. We now consider a new solution $\bX_{\alpha}=\bX+\alpha\mb{D}$, where $\displaystyle 0<\alpha\leq\frac{x_{ij}}{MNr}$. Clearly, $\bX_{\alpha}\geq\mb{0}$. Thus we have: $\norm{\bX_{\alpha}}_1=\norm{\bX}_1+\alpha MN(1-r)$. Applying the triangle inequality, we can bound $\norm{\bX_{\alpha}}_*$ as follows:
$$
\norm{\bX_{\alpha}}_*\leq\norm{\bX}_*+\alpha\norm{\mb{D}}_*\leq\norm{\bX}_*+\alpha(\sqrt{MN}+MNr).
$$
Thus we have:
$$
\norm{\bX_{\alpha}}_\theta\leq\norm{\bX}_{\theta}+\alpha\sqrt{MN}\left[(1+r\sqrt{MN})+\theta\sqrt{MN}(1-r)\right].
$$
Since $\displaystyle\theta>\theta_B=\frac{1}{\sqrt{MN}}\left(\frac{\bar{a}({\cal I},{\cal J})\sqrt{MM}+a_{\max}(\overline{{\cal I},{\cal J}})}{\bar{a}({\cal I},{\cal J})-a_{\max}(\overline{{\cal I},{\cal J}})}\right)$ and $0<a_{ij}\leq a_{\max}(\overline{{\cal I},{\cal J}})<\bar{a}({\cal I},{\cal J})$, we have:
$$
\theta>\frac{1}{\sqrt{MN}}\left(\frac{1+r\sqrt{MN}}{r-1}\right).
$$
This implies that $\norm{\bX_{\alpha}}_{\theta}<\norm{\bX}_{\theta}=1$ for all $\displaystyle 0<\alpha\leq\frac{x_{ij}}{MNr}$. Now consider the scaled solution $\displaystyle\bX_{\alpha}^s=\frac{1}{\norm{\bX_{\alpha}}_{\theta}}\bX_{\alpha}$, which is also a feasible solution of Problem \refs{eq:dnorm}. In terms of the objective, we have: $\inprod{\mb{A}}{\mb{D}}=\mb{e}_{\cal I}^T\mb{A}\mb{e}_{\cal J}-MNra_{ij}=0$. Thus $\inprod{\mb{A}}{\bX_{\alpha}}=\inprod{\mb{A}}{\bX}=\norm{\mb{A}}_{\theta}^*$ for all $\alpha$. We then have: 
$$
\inprod{\mb{A}}{\bX_{\alpha}^s}=\frac{\norm{\mb{A}}_{\theta}^*}{\norm{\bX_{\alpha}}_{\theta}}>\norm{\mb{A}}_{\theta}^*,\quad\forall\,\alpha\in\left(0,\frac{x_{ij}}{MNr}\right],
$$
which is a contradiction because $\norm{\mb{A}}_{\theta}^*$ is the optimal value of Problem \refs{eq:dnorm}. 

Thus we can conclude that if $\mb{A}\geq\mb{0}$, all optimal solutions $\bX$ of Problem \refs{eq:dnorm} are sparse with $x_{ij}=0$ for all $(i,j)\not\in({\cal I},{\cal J})$ when $\theta>\theta_B$. 
\end{pf}


\section{Random noise}
\label{sec:randomnoise}

The main technical result of this article is that the proposed algorithm
can find a large rank-one submatrix hidden in a substantial amount of
noise.  The noise takes two forms: the rank-one submatrix itself has 
random noise
added to it (so that its rank is no longer 1), and the entries outside
the rank-one submatrix are generated by a random process.

First, we recall the following definition: a random variable
$x$ is {\em $b$-subgaussian} if its mean is zero, and if there
exists a $b>0$ such that for all $t\ge 0$,
\begin{equation}
\mathbb{P}(|x|\ge t)\le \exp(-t^2/(2b^2)).
\label{eq:subgau}
\end{equation}
For example, a normally distributed variable or any mean-zero
variable with a discrete distribution is subgaussian.

The result of this section is the following bound.
For this entire section, we adopt the notation that $\mb{e}_M$ denotes
the vector of all 1's of length $M$, and similarly for $\mb{e}_N$.
\begin{theorem}
Let $\mb{A}$ be an $m\times n$ matrix defined as follows.
\begin{equation}
\mb{A} = \left(
\ba{cc}
\sigma_0 \mb{u}_0\mb{v}_0^T & \mb{0} \\
\mb{0} & \mb{0}
\ea
\right) +
\left(
\ba{cc}
\mb{R}_{11} & \mb{R}_{12} \\ 
\mb{R}_{21} & \mb{R}_{22}
\ea
\right),
\label{eq:Aform}
\end{equation}
where $\sigma_0>0$, $\mb{u}_0\in\R^M$,  
$\mb{u}_0\ge\mb{0}$,
$M<m$, and
$\mb{v}_0\in\R^N$, 
$\mb{v}_0\ge\mb{0}$,
$N<n$.  Furthermore, assume that
$\mb{u}_0=\mb{e}_M+\mb{p}$ with $\norm{\mb{p}}_2\leq c_1\sqrt{M}$, 
and 
$\mb{v}_0=\mb{e}_N+\mb{q}$ with
$\norm{\mb{q}}_2\leq c_2\sqrt{N}$.
The matrix
$\mb{R}$ is a random matrix with i.i.d.\ nonnegative elements $r_{ij}$
with mean $c_3\sigma_0$, where $c_3>0$ is a constant, such that
$r_{ij}/\sigma_0-c_3$ is $b$-subgaussian.
Here $c_1,c_2,c_3$ are positive constants.
Assume that these scalar constants
$c_1,c_2,c_3$ satisfy the following relations
\begin{equation}
c_5 \le 1/3,\qquad c_3+c_5<1.
\label{eq:cond1}
\end{equation}
where $c_5$ is chosen to satisfy
\begin{equation}
c_5>c_1+c_2+c_1c_2.
\label{eq:c5def}
\end{equation}

Under these hypotheses concerning $\mb{A}$, 
and assuming $\theta$ satisfies
\begin{eqnarray}
\theta &\le& \min\left(\frac{1}{c_3+c_5}, \frac{1+c_3-3c_5}{2c_5}\right)
\cdot \frac{1}{\sqrt{MN}},
\label{eq:thetacond1}
\\
\theta &\ge &\frac{2c_3}{1-c_3-c_5}\cdot \frac{1}{\sqrt{MN}},
\label{eq:thetacond2}
\end{eqnarray}
the solution $\mb{X}$ to problem \req{eq:LAROSrelaxation} 
is a rank-one
matrix with positive entries in positions that are indexed by $\{1,\ldots,M\}\times\{1,\ldots,N\}$ and zeros elsewhere with probability exponentially
close to $1$ (i.e., of the form $1-\exp(-(M+N)^{\rm const})$) provided
that $MN\ge \Omega\left((M+N)^{4/3}\right)$ and
$MN\ge \Omega(m+n)$.  Here, the constants implicit in the $\Omega(\cdot)$
notation depend on $b$ and $c_5$.  See \req{eq:MN1}--\req{eq:k2def} 
below for a detailed presentation of these constants.
\end{theorem}

\noindent
{\bf Remarks.}  

\begin{enumerate}
\item
Naturally, the theorem also applies if the $MN$ distinguished
entries occur as any $M\times N$ submatrix of $\mb{A}$; we have numbered the
distinguished submatrix first in order to simplify notation.

\item
It is not enough to assume simply that $\mb{u}_0>\mb{0}$
and $\mb{v}_0>\mb{0}$ because if these vectors have very small entries,
then they cannot be distinguished from the noise.

\item
This theorem is not a consequence of Theorem~\ref{thm:zblock}
because the hypotheses do not force entries outside the
distinguished block to be smaller than the average of the 
distinguished block's entries.

\item
The relationships among the constants as well as \req{eq:thetacond1},
\req{eq:thetacond2} can all be satisfied provided $c_3,c_5$ are sufficiently
small.

\item
The result holds with probability exponentially close
to 1 as long as
$M\sim N$ and $M\ge \Omega(m^{1/2})$, $N\ge \Omega(n^{1/2})$.  Thus, the 
rank-one submatrix can be much smaller than the entire matrix $\mb{A}$.
\end{enumerate}

Before beginning the proof of the theorem, we require the following
key lemma regarding 
matrices constructed from independent $b$-subgaussian random variables.
\begin{lemma}
\label{lem:2norm}
Let $\mb{B}\in\R^{m\times n}$ be a random matrix, where $b_{ij}$ are independent $b$-subgaussian random variables for all $i=1,\ldots,m$, and $j=1,\ldots,n$. 
Then for any $u>0$,
\ben
\item[(i)] $\displaystyle\mathbb{P}\left(\norm{\mb{B}}\geq u\right)\leq \exp\left(-\left(\frac{8u^2}{81b^2}-(\log 7)(m+n)\right)\right)$
\item[(ii)]
 $\displaystyle\mathbb{P}\left(\norm{\mb{C}\mb{B}}\geq u\right)\leq
\exp\left(-\left(\frac{8u^2}{81b^2\norm{\vect{C}}^2}-(\log 7)
(m+n)\right)\right)$, where $\mb{C}$ is a deterministic matrix.
\een
\end{lemma}
The proof of this lemma follows the proof techniques by Litvak 
et al.\ \cite{Litvak05}. Major steps are shown as follows.

\begin{pf}
\ben
\item[(i)] We have: $\displaystyle\norm{\mb{B}}=\max_{\norm{\bx}_2=\norm{\by}_2=1}\by^T\mb{B}\bx$. We discretize the unit balls in $\R^n$ and $\R^m$ by finite $\eps$-nets, where $\eps\in (0,1)$. An $\eps$-net of a set $\cal K$ is the subset $\cal N$ such that for all $\bx\in\cal K$, there exists $\by\in\cal N$ such that $\norm{\bx-\by}_2\leq \eps$. Using a construction proof, we can prove that there exists a finite $\eps$-net of the unit ball in $\R^n$ with the cardinality of no more than $\displaystyle\left(\frac{2}{\eps}+1\right)^n$. Let $\cal N$ and $\cal M$ be the finite $\eps$-nets of the unit balls in $\R^n$ and $\R^m$ with minimum cardinality, respectively. Applying the triangle inequality, we have:
$$
\norm{\mb{B}}\leq\frac{1}{(1-\eps)^2}\max_{\bx\in{\cal N},\by\in{\cal M}}\by^T\mb{B}\bx.
$$
We can bound the tail probability $P(\norm{\mb{B}}\geq u)$ as follows.
$$
\mathbb{P}\left(\norm{\mb{B}}\geq u\right)\leq\left(\frac{2}{\eps}+1\right)^{m+n}\max_{\bx\in{\cal N},\by\in{\cal M}}\mathbb{P}\left(\by^T\mb{B}\bx\geq (1-\eps)^2u\right).
$$

We have, $b_{ij}$ are independent $b$-subgaussian random variables; therefore, $\displaystyle\by^T\mb{B}\bx=\sum_{i=1}^m\sum_{j=1}^n\left(x_jy_i\right)b_{ij}$ is also a $b$-subgaussian random variable since $\displaystyle \sum_{i=1}^n\sum_{j=1}^n\left(x_jy_i\right)^2=\norm{\bx}_2^2\norm{\by}_2^2=1$. Thus we have:
$$
\mathbb{P}\left(\norm{\mb{B}}\geq u\right)\leq \left(\frac{2}{\eps}+1\right)^{m+n}e^{-\frac{(1-\eps)^4u^2}{2b^2}}.
$$
Letting $\eps=1/3$, we obtain the inequality
$$
\mathbb{P}\left(\norm{\mb{B}}\geq u\right)\leq e^{-\left(\frac{8u^2}{81b^2}-(\log 7)(m+n)\right)}
$$
\item[(ii)] We have: $\by^T\mb{C}\mb{B}\bx=\left(\mb{C}^T\by\right)^T\mb{B}\bx$, thus:
$$
\mathbb{P}\left(\by^T\mb{C}\mb{B}\bx\geq u\right)\leq e^{-\frac{u^2}{2b^2\norm{\vect{C}^T\by}_2^2\norm{\bx}_2^2}}\leq e^{-\frac{u^2}{2b^2\norm{\vect{C}}^2}}.
$$
Applying similar arguments, we can then obtain the inequality (ii)
of the Lemma.
\een
\end{pf}

We now turn to the proof of the main theorem.
We now would like to find conditions on $\theta$ and the constants 
so that Problem \refs{eq:LAROSrelaxation} has an optimal solution $\bX$ of the form
$$
\bX = \left(
\ba{cc}
\sigma_1\mb{u}_1\mb{v}_1^T & \mb{0}\\
\mb{0} & \mb{0}
\ea
\right),
$$
where $\mb{u}_1>\mb{0}$, $\norm{\mb{u}_1}_2=1$, and $\mb{v}_1>\mb{0}$,
$\norm{\mb{v}_1}_2=1$. If $\mb{u}_1$ and $\mb{v}_1$ are determined,
$\sigma_1$ can be easily calculated in order to satisfy the condition
$\inprod{\mb{A}}{\bX}=1$ of the optimal solution.  Thus the main task is
to find $\mb{u}_1$ and $\mb{v}_1$ if they exist. We will construct
them using optimality conditions derived in the previous section for
Problem \refs{eq:dnorm} (equivalent to Problem \refs{eq:LAROSrelaxation}) and
its dual, Problem \refs{eq:gfunc}. Defining
$\mb{u}=[\mb{u}_1;\mb{0}]\in\R^m$ and
$\mb{v}=[\mb{v}_1;\mb{0}]\in\R^n$, we can then write the optimality
conditions as follows:
\begin{quote}
There exists $\bY$ and $\bZ$ such that $\mb{A}=\bY+\bZ$ and
$$
\bY = \norm{\mb{A}}_{\theta}^*(\mb{u}\mb{v}^T+\mb{W}),\quad\bZ = \theta\norm{\mb{A}}_{\theta}^*\mb{V},
$$
where $\norm{\mb{W}}\leq 1$, $\mb{W}^T\mb{u}=\mb{0}$, $\mb{W}\mb{v}=\mb{0}$, and $\norm{\mb{V}}_{\infty}\leq 1$, $\mb{V}_{11}=\mb{e}_M\mb{e}_N^T$.
\end{quote}
These conditions come from the properties of the subgradient $\partial\norm{\bY}$ and $\partial\norm{\bZ}_{\infty}$ and the fact that $\bX$ belongs to these sets (up to appropriate scaling factors).

In the following analysis, we will construct 
$(\mb{V},\mb{W})$ so that the optimality conditions are
satisfied.  The entries of these matrices will be constructed
separately for the four subblocks of $\mb{A}$, starting with
the $(1,1)$ block.  Breaking the equation $\mb{Y}+\mb{Z}=\mb{A}$
into blocks and scaling by $1/\Vert \mb{A}\Vert_\theta^*$, we obtain
the following more detailed optimality conditions.
\begin{eqnarray}
\mb{u}_1\mb{v}_1^T + \mb{W}_{11}+\theta\mb{e}_M\mb{e}_N^T&=&\mb{A}_{11}/\Vert \mb{A}\Vert_\theta^*=
(\sigma_0\mb{u}_0\mb{v}_0^T+\mb{R}_{11})/\Vert \mb{A}\Vert_\theta^*, 
\label{eq:vw11}\\
\mb{W}_{ij}+\theta \mb{V}_{ij}&=&\mb{A}_{ij}/\Vert \mb{A}\Vert_\theta^*=\mb{R}_{ij}/\Vert \mb{A}\Vert_\theta^*,\label{eq:vwij}
\end{eqnarray}
where the second line applies to $(i,j)$ equal to $(1,2)$, $(2,1)$ 
and $(2,2)$.
Following this block matrix notation, the
remaining optimality conditions to be established are
$\mb{W}_{11}^T\mb{u}_1=\mb{0}$,
$\mb{W}_{21}^T\mb{u}_1=\mb{0}$,
$\mb{W}_{11}\mb{v}_1=\mb{0}$,
$\mb{W}_{12}\mb{v}_1=\mb{0}$,
$\Vert\mb{V}\Vert_\infty\le 1$,
$\Vert\mb{W}\Vert\le 1$.
We shall establish the latter inequality by proving more
specifically that 
$\Vert\mb{W}_{ij}\Vert\le 1/2$ for $(i,j)\in\{1,2\}\times\{1,2\}$.

We begin with the $(1,1)$ block of this equation.  The conditions
$\norm{\mb{W}_{11}}\leq 1/2$, 
$\mb{W}_{11}^T\mb{u}_1=\mb{0}$, $\mb{W}_{11}\mb{v}_1=\mb{0}$
imply that
$\norm{\mb{u}_1\mb{v}_1^T+\mb{W}_{11}}=1$.  This is because
the dominant singular triple of $\mb{u}_1\mb{v}_1^T+\mb{W}_{11}$
must be $(1,\mb{u}_1,\mb{v}_1)$ by the conditions.
Equivalent to $\norm{\mb{u}_1\mb{v}_1^T+\mb{W}_{11}}=1$ is
$$
\norm{\frac{1}{\norm{A}_{\theta}^*}\mb{A}_{11}-\theta\mb{V}_{11}}=1,
$$
where, as noted earlier, we are required to take $\mb{V}_{11}=\mb{e}_M\mb{e}_N^T$.

Thus the first necessary condition for $\mb{X}$ to be the optimal solution
is that there exists $\lambda>0$ such that
$f(\lambda)=\norm{\lambda\mb{A}_{11}-\theta\mb{V}_{11}}=1$. If
such a $\lambda$ is identified, then $\mb{u}_1$ and $\mb{v}_1$
can be easily found since $\mb{u}_1\mb{v}_1^T$ is the rank-one
approximation of $\lambda\mb{A}_{11}-\theta\mb{V}_{11}$. Note that the
nonnegativity of $\mb{u}_1$ and $\mb{v}_1$ will require additional
conditions which will be discussed later. We have:
\begin{align*}
\lambda\mb{A}_{11}-\theta\mb{V}_{11}&= \lambda\left[\sigma_0\mb{u}_0\mb{v}_0^T+\mb{R}_{11}\right]-\theta\mb{e}_M\mb{e}_N^T\\
\quad &= \lambda\left[\sigma_0(\mb{e}_M+\mb{p})(\mb{e}_N+\mb{q})^T+\mb{R}_{11}\right]-\theta\mb{e}_M\mb{e}_N^T\\
\quad &=(\lambda\sigma_0-\theta)\mb{e}_M\mb{e}_N^T+\lambda\left[\sigma_0(\mb{e}_M\mb{q}^T+\mb{p}\mb{e}_N^T+\mb{p}\mb{q}^T)+\mb{R}_{11}\right]\\
\quad &=[\lambda\sigma_0(1+c_3)-\theta]\mb{e}_M\mb{e}_N^T+\lambda\left[\sigma_0(\mb{e}_M\mb{q}^T+\mb{p}\mb{e}_N^T+\mb{p}\mb{q}^T)+(\mb{R}_{11}-c_3\sigma_0\mb{e}_M\mb{e}_N^T)\right].
\end{align*}

We have: $f(\lambda)\rightarrow+\infty$ when $\lambda\rightarrow+\infty$ since $\mb{A}_{11}\ne \mb{0}$. 
Now define $\displaystyle\lambda_0=\frac{\theta}{\sigma_0(1+c_3)}$ to
make the first term vanish, yielding
\begin{align*}
\lambda_0\mb{A}_{11}-\theta\mb{V}_{11} &=\lambda_0\left[\sigma_0(\mb{e}_M\mb{q}^T+\mb{p}\mb{e}_N^T+\mb{p}\mb{q}^T)+(\mb{R}_{11}-c_3\sigma_0\mb{e}_M\mb{e}_N^T)\right]\\
\quad & = \lambda_0\left[\sigma_0\mb{P}+\mb{Q}\right],
\end{align*}
where $\mb{P}=\mb{e}_M\mb{q}^T+\mb{p}\mb{e}_N^T+\mb{p}\mb{q}^T$ and $\mb{Q}=\mb{R}_{11}-c_3\sigma_0\mb{e}_M\mb{e}_N^T$. We now bound the spectral norm of $\mb{P}$ and $\mb{Q}$ as follows.
\begin{align*}
\norm{\mb{P}}&\leq\norm{\mb{e}_M\mb{q}^T}+\norm{\mb{p}\mb{e}_N^T}+\norm{\mb{p}\mb{q}^T}\\
\quad & =\norm{\mb{e}}_2\norm{\mb{q}}_2+\norm{\mb{p}}_2\norm{\mb{e}_N}_2+\norm{\mb{p}}_2\norm{\mb{q}}_2\\
\quad & \leq c_1\sqrt{MN}+c_2\sqrt{MN}+c_1c_2\sqrt{MN}\\
\quad & = (c_1+c_2+c_1c_2)\sqrt{MN}.
\end{align*}

Matrix $\mb{Q}/\sigma_0$ is random with i.i.d.\ elements that are $b$-subgaussian. Thus by Lemma~\ref{lem:2norm}(i),
$$
\mathbb{P}\left(\norm{\mb{Q}}\geq u\sigma_0\right)\leq \exp\left(-\left(\frac{8u^2}{81b^2}-(\log 7)(M+N)\right)\right),
$$
for any $u>0$.  Let us fix $u=(MN)^{3/8}$ to obtain
\begin{equation}
\mathbb{P}\left(\norm{\mb{Q}}\geq \sigma_0(MN)^{3/8}\right)\leq \exp\left(-\left(\frac{(MN)^{3/4}}{81b^2}-(\log 7)(M+N)\right)\right).
\label{eq:probq}
\end{equation}
For the remainder of this analysis, we will impose
the assumption that the event
in \req{eq:probq} does not happen.  At the end of the proof the
right-hand side
\req{eq:probq} will be one of the terms in the failure probability
of identifying the optimal $\mb{X}$.

Thus, $\Vert\mb{Q}\Vert\le o(1)\sigma_0\sqrt{MN}$, so
applying the triangle inequality,
\begin{eqnarray}
\Vert\sigma_0 \mb{P}+\mb{Q}\Vert &\le& \sigma_0(c_1+c_2+c_1c_2+o(1))\sqrt{MN}\nonumber\\
&\le& \sigma_0c_5\sqrt{MN}.\label{eq:sigma0PQ}
\end{eqnarray}
by \req{eq:c5def}. (The strict inequality `$>$' in \req{eq:c5def} is used
in order to absorb the $o(1)$ term.)
Therefore,
\begin{eqnarray*}
f(\lambda_0)&=& \lambda_0\Vert \sigma_0\mb{P}+\mb{Q}\Vert \\
&\le&\frac{c_5\theta\sqrt{MN}}{1+c_3}.
\end{eqnarray*}
Thus with high probability, $f(\lambda_0)\leq 1$ if 
\begin{equation}
\theta < \left(\frac{1+c_3}{c_5}\right)\frac{1}{\sqrt{MN}},
\label{eq:thetaub1a}
\end{equation}
Inequality \req{eq:thetaub1a} is a consequence of 
\req{eq:thetacond1} stated in the theorem.
This inequality implies
$f(\lambda_0)\leq 1$, and, due to the continuity of function $f$,
there exists $\lambda^*\geq\lambda_0$ such that
$f(\lambda^*)=1$. We will prove that under some additional conditions,
this value $\lambda^*$ satisfies all other optimality conditions of
Problem \refs{eq:LAROSrelaxation} and indeed
$\displaystyle\norm{\mb{A}}_{\theta}^*=\frac{1}{\lambda^*}$.

Let us recall that $\Vert\lambda^*\mb{A}_{11}-\theta \mb{V}_{11}\Vert=1$,
i.e., $\Vert(\lambda^*\sigma_0(1+c_3)-\theta)\mb{e}_M\mb{e}_N^T+\lambda^*(\sigma_0
\mb{P}+\mb{Q})\Vert =1$.  Applying the fact that 
$\Vert \mb{e}_M\mb{e}_N^T\Vert=\sqrt{MN}$
and the triangle inequality 
twice
to this equation yields
$$
[\lambda^*\sigma_0(1+c_3)-\theta]\sqrt{MN} - \lambda^*\norm{\sigma_0\mb{P}+\mb{Q}}\leq 1\leq[\lambda^*\sigma_0(1+c_3)-\theta]\sqrt{MN} + \lambda^*\norm{\sigma_0\mb{P}+\mb{Q}}.
$$ 
Applying \req{eq:sigma0PQ} yields
$$
[\lambda^*\sigma_0(1+c_3-c_5)-\theta]\sqrt{MN}\leq 1\leq [\lambda^*\sigma_0(1+c_3+c_5)-\theta]\sqrt{MN}.
$$
Rearranging this chain of inequalities and using the fact that
$1+c_3-c_5>0$, which follows
from \req{eq:cond1} stated in the theorem, yields
\begin{equation}
\frac{1+\theta\sqrt{MN}}{\sigma_0(1+c_3+c_5)\sqrt{MN}}\leq\lambda^*\leq\frac{1+\theta\sqrt{MN}}{\sigma_0(1+c_3-c_5)\sqrt{MN}},
\label{eq:lambdastarbounds}
\end{equation}
with high probability.

We wish to establish that $\lambda^*\sigma_0-\theta\ge 0$.  Using the left 
inequality in \req{eq:lambdastarbounds} yields:
\begin{eqnarray*}
\lambda^*\sigma_0-\theta
& \ge &
\frac{1+\theta\sqrt{MN}}{(1+c_3+c_5)\sqrt{MN}}-\theta
\\
& = &
\frac{1-(c_3+c_5)\theta\sqrt{MN}}{(1+c_3+c_5)\sqrt{MN}}.
\end{eqnarray*}
Thus, nonnegativity of $\lambda^*\sigma_0-\theta$ is implied by
the inequality $\theta \le 1/((c_3+c_5)\sqrt{MN})$, which is
a consequence of assumption \req{eq:thetacond1}.

Since $\lambda^*\sigma_0-\theta\geq 0$,
$$
\lambda^*\mb{A}_{11}-\theta\mb{V}_{11}=(\lambda^*\sigma_0-\theta)\mb{e}_M\mb{e}_N^T+\lambda^*(\sigma_0\mb{P}+\mb{R}_{11})>\mb{0}.
$$
Applying the Perron-Frobenius theorem, we obtain the positivity of $\mb{u}_1$ and $\mb{v}_1$.

We also need
$\norm{\mb{W}_{11}}\le 1/2$.  Recall
$\norm{\mb{W}_{11}}=\sigma_2(\lambda^*\mb{A}_{11}-\theta\mb{V}_{11})$,
the second largest singular value of
$\lambda^*\mb{A}_{11}-\theta\mb{V}_{11}$, since 
$\lambda^*\mb{A}_{11}-\theta\mb{V}_{11}=\mb{u}_1\mb{v}_1^T+\mb{W}_{11}$.
Using the well-known fact
that 
$$\sigma_2(\mb{A})=
\min\{\norm{\mb{A}-\mb{S}}:\mathop{\rm rank}(\mb{S})\le 1\},$$
we obtain
$$
\norm{\mb{W}_{11}}\leq\norm{\lambda^*(\sigma_0\mb{P}+\mb{Q})}.
$$
Here we selected $\mb{S}$ to be $[\lambda^*\sigma_0(1+c_3)-\theta]\mb{e}_M\mb{e}_N^T$. With high probability, we obtain the bound 
$\norm{\mb{W}_{11}}\leq\lambda^*\sigma_0c_5\sqrt{MN}$
from \req{eq:sigma0PQ}.

Using the upper bound on $\lambda^*$ from \req{eq:lambdastarbounds}, we have:
$$
\norm{\mb{W}_{11}}\leq \frac{(1+\theta\sqrt{MN})c_5}{1+c_3-c_5}.
$$

In order to obtain $\norm{\mb{W}_{11}}\leq 1/2$, a sufficient condition is
\begin{equation}
\frac{(1+\theta\sqrt{MN})c_5}{1+c_3-c_5}\leq\frac{1}{2}
\label{eq:thetac}
\end{equation}
which is rearranged as
$$
\theta\leq\left(\frac{1+c_3-3c_5}{2c_5}\right)\frac{1}{\sqrt{MN}}.
$$
The latter inequality follows from \req{eq:thetacond1};
the numerator of the right-hand side is positive by
\req{eq:cond1}.

Turning to \req{eq:vwij} when $(i,j)=(2,2)$,
we need to find $\mb{W}_{22}$ and $\mb{V}_{22}$ that satisfy
$$
\lambda^*\mb{R}_{22}=\mb{W}_{22}+\theta\mb{V}_{22},
$$
$\norm{\mb{W}_{22}}\le 1/2$,
and $\norm{\mb{V}_{22}}_{\infty}\leq 1$. Consider the assignment
$\displaystyle\mb{V}_{22}=\frac{\lambda^*\sigma_0c_3}{\theta}\mb{e}_{m-M}\mb{e}_{n-N}^T$
and $\mb{W}_{22}=\lambda^*\mb{R}_{22}-\theta\mb{V}_{22}$.  The
coefficient $\lambda^*\sigma_0c_3/\theta$ is chosen for the
definition of $\mb{V}_{22}$  so that the
entries of the remainder term $\mb{W}_{22}$
have mean zero.

 The requirement $\norm{\mb{V}_{22}}_{\infty}\leq 1$ is satisfied if and only if ${\lambda^*\sigma_0c_3}/{\theta}\leq 1$. Because of the upper bound on
$\lambda^*$ established by
\req{eq:lambdastarbounds},
this requirement is satisfied if
$$
\frac{(1+\theta\sqrt{MN})c_3}{\theta(1+c_3-c_5)\sqrt{MN}}\leq 1\Leftrightarrow\theta\geq\left(\frac{c_3}{1-c_5}\right)\frac{1}{\sqrt{MN}},\quad c_5<1.
$$
This inequality is assured by \req{eq:cond1} and \req{eq:thetacond2}.
(In particular, \req{eq:cond1} implies $c_5<1$.)

To bound $\Vert \mb{W}_{22}\Vert$, 
consider $\mb{W}_{22}/(\lambda^*\sigma_0)$, 
which is a random matrix with i.i.d.\ elements that are $b$-subgaussian. 
Applying Lemma~\ref{lem:2norm}(i)
to $\mb{W}_{22}/(\lambda^*\sigma_0)$
and taking $u=1/(2\lambda^*\sigma_0)$ yields
$$
\mathbb{P}\left(\norm{\mb{W}_{22}}\geq 1/2\right)\leq \exp\left(-\left(\frac{2}{81b^2(\lambda^*\sigma_0)^2}-(\log 7)(m-M+n-N)\right)\right).
$$ 
Use the upper bound on $\lambda^*$ from 
\req{eq:lambdastarbounds} to obtain
the following tail bound:
$$
\mathbb{P}\left(\norm{\mb{W}_{22}}\geq \frac{1}{2}\right)\leq \exp\left(-\left(\frac{2(1+c_3-c_5)^2}{81b^2(1+\theta\sqrt{MN})^2}MN-(\log 7)(m-M+n-N)\right)\right).
$$
{From} \req{eq:thetac} we obtain
\begin{equation}
\frac{(1+c_3-c_5)^2}{(1+\theta\sqrt{MN})^2}\ge 4c_5^2
\label{eq:sqtheta}
\end{equation}
hence
\begin{equation}
\mathbb{P}\left(\norm{\mb{W}_{22}}\geq \frac{1}{2}\right)\leq \exp\left(-\left(\frac{8c_5^2MN}{81b^2}-(\log 7)(m-M+n-N)\right)\right).
\label{eq:probw22}
\end{equation}

Now consider \req{eq:vwij}
when $(i,j)=(1,2)$.
Again we need to find $\mb{W}_{12}$ and $\mb{V}_{12}$ such that
$$
\lambda^*\mb{R}_{12}=\mb{W}_{12}+\theta\mb{V}_{12},
$$
$\norm{\mb{W}_{12}}\le 1/2$, 
$\norm{\mb{V}_{12}}_{\infty}\leq 1$, and
$\mb{W}_{12}^T\mb{u}_1=\mb{0}$. We construct $\mb{W}_{12}$ and
$\mb{V}_{12}$ column by column as follows:
\begin{eqnarray*}
\mb{V}_{12}(:,i)&=&\frac{\lambda^*\mb{R}_{12}(:,i)^T\mb{u}_1}{\theta\norm{\mb{u}_1}_1}\mb{e}_{M},\quad
\mb{W}_{12}(:,i)=\lambda^*\mb{R}_{12}(:,i)-\theta\mb{V}_{12}(:,i)
\end{eqnarray*}
all $i=1,\ldots,n-N$. By construction we have $\mb{W}_{12}(:,i)^T\mb{u}_1=0$
for all $i=1,\ldots,n-N$. Now consider the requirement that
$\norm{\mb{V}_{12}(:,i)}_{\infty}\leq 1$ for all
$i=1,\ldots,n-N$. The requirement is equivalent to
$$\frac{\lambda^*\mb{R}_{12}(:,i)^T\mb{u}_1}{\theta\norm{\mb{u}_1}_1}\leq
1$$ for all $i=1,\ldots,n-N$.
Subtract $\lambda^*c_3\sigma_0/\theta$ from both sides and apply the
identity $\mb{e}_M^T\mb{u}_1=\Vert\mb{u}_1\Vert_1$ to obtain
$$\frac{\lambda^*(\mb{R}_{12}(:,i)^T-c_3\sigma_0\mb{e}_M^T)\mb{u}_1}{\theta\norm{\mb{u}_1}_1}\leq
1-\frac{\lambda^*c_3\sigma_0}{\theta}.$$

We will establish this inequality in two steps.  First,
we establish that $\lambda^*c_3\sigma_0/\theta\le 1/2$.  Because
of \req{eq:lambdastarbounds}, it suffices to establish that
\begin{equation}
\frac{c_3(1+\theta\sqrt{MN})}{\theta(1+c_3-c_5)\sqrt{MN}}\le \frac{1}{2}.
\label{eq:c3theta12}
\end{equation}
This can be rearranged into
$$\theta \ge \frac{2c_3}{(1-c_3-c_5)\sqrt{MN}},$$
which follows from \req{eq:thetacond2}.

Second, we establish that with probability exponentially close to 1,
$$\frac{\lambda^*(\mb{R}_{12}(:,i)^T-c_3\sigma_0\mb{e}_M^T)\mb{u}_1}{\theta\norm{\mb{u}_1}_1}\le \frac{1}{2}.$$
Notice that $r_{ji}/\sigma_0-c_3$ is $b$-subgaussian; thus, 
$$\frac{1}{\sigma_0}\left(\mb{R}_{12}(:,i)-c_3\sigma_0\mb{e}_M\right)^T\mb{u}_1
=\frac{1}{\sigma_0}\sum_{j=1}^Mu_1(j)\left(R_{12}(j,i)-c_3\sigma_0\right)$$
is also $b$-subgaussian since $\norm{\mb{u}_1}_2=1$. 
Thus, by \req{eq:subgau},
taking $x=
\left(\mb{R}_{12}(:,i)-c_3\sigma_0\mb{e}_M\right)^T\mb{u}_1/\sigma_0$ and
taking $t=\theta\Vert\mb{u}_1\Vert_1/(2\lambda^*\sigma_0)$,
\begin{eqnarray}
\mathbb{P}\left(\frac{\lambda^*(\mb{R}_{12}(:,i)^T-c_3\sigma_0\mb{e}_M)\mb{u}_1}{\theta\norm{\mb{u}_1}_1}> \frac{1}{2}\right)
&\le&
\exp\left(-\theta^2\Vert\mb{u}_1\Vert_1^2/(8b^2(\lambda^*\sigma_0)^2)\right). 
\nonumber\\
&\le&
\exp\left(-c_3^2\Vert\mb{u}_1\Vert_1^2/(2b^2)\right).
\label{eq:plr-c}
\end{eqnarray}
since, as noted above $\theta/(\lambda^*\sigma_0)\ge 2c_3$.

To proceed, we now need a lower bound for $\norm{\mb{u}_1}_1$.
Let $\mb{F}$ denote
$\lambda^*\mb{A}_{11}-\theta\mb{V}_{11}$, which is equal to
$[\lambda^*\sigma_0(1+c_3)-\theta]\mb{e}_M\mb{e}_N^T+\lambda^*(\sigma_0\mb{P}+\mb{Q}).$
We
know that $\mb{u}_1$ is the first (left) singular vector of
$\mb{F}$.  Letting
$\mb{X}_0=[\lambda^*\sigma_0(1+c_3)-\theta]\mb{e}_M\mb{e}_N^T$
and
$\mb{E}=\lambda^*(\sigma_0\mb{P}+\mb{Q})$, we
then have $\mb{F}=\mb{X}_0+\mb{E}$, and
$\mb{X}_0$ is a rank-one matrix with a single nonzero
singular value equal to $(\lambda^*\sigma_0(1+c_3)-\theta)\sqrt{MN}$
and with left singular vector $\mb{e}_M/\sqrt{M}$ and right
singular vector $\mb{e}_N/\sqrt{N}$. 
Furthermore, since $\Vert \mb{F}\Vert=1$, we know that the singular
value of $\mb{X_0}$ is at least $1-\Vert\mb{E}\Vert$ by Corollary
8.6.2 of \cite{GVL}.
Thus, by Theorem 8.6.5 of \cite{GVL}, 
\begin{equation}
\left\Vert \mb{u}_1-\frac{\mb{e}_M}{\sqrt{M}}\right\Vert
\le 
\frac{4\Vert \mb{E}\Vert}{1-\Vert \mb{E}\Vert}\le 4/5,
\label{eq:emu1}
\end{equation}
provided that $\Vert \mb{E}\Vert\le 1/6$.  Thus, the next step in
the analysis is to show that $\Vert \mb{E}\Vert\le 1/6$.  This follows
from the following sequence of inequalities:
\begin{eqnarray*}
\Vert \mb{E}\Vert &=&\lambda^*\Vert \sigma_0\mb{P}+\mb{Q}\Vert \\ 
&\le & 
\frac{c_5(1+\theta\sqrt{MN})}{1+c_3-c_5} \\
&\le &
1/6,
\end{eqnarray*}
where the second line holds with high probability according to
\req{eq:sigma0PQ} and the third line follows
from \req{eq:thetac} and \req{eq:cond1}.

Thus, we have established that $\Vert\mb{E}\Vert\le 1/6$
with high probability, which in
turn implies that
\begin{eqnarray*}
\Vert\mb{u}_1\Vert_1 & = &
\mb{e}_M^T\mb{u}_1 \\
& = &
(\mb{e}_M-M^{1/2}\mb{u}_1+M^{1/2}\mb{u}_1)^T\mb{u}_1 \\
& \ge &
M^{1/2}\mb{u}_1^T\mb{u}_1 - |(\mb{e}_M-M^{1/2}\mb{u}_1)^T\mb{u}_1| \\
& \ge &
M^{1/2}\mb{u}_1^T\mb{u}_1  - \Vert\mb{e}_M-M^{1/2}\mb{u}_1\Vert\cdot\Vert \mb{u}_1\Vert \\
& = &
M^{1/2} - M^{1/2}\Vert M^{-1/2}\mb{e}_M-\mb{u}_1\Vert \\
&\ge&
M^{1/2}-(4/5)M^{1/2},
\end{eqnarray*}
where the last line is obtained from \req{eq:emu1}.
This gives a lower bound of $M^{1/2}/5$ on $\Vert \mb{u}\Vert_1$.
Thus, substituting this into 
\req{eq:plr-c}
yields
$$
\mathbb{P}\left(\frac{\lambda^*(\mb{R}_{12}(:,i)^T-c_3\sigma_0\mb{e}_M)\mb{u}_1}{\theta\norm{\mb{u}_1}_1}> \frac{1}{2}\right)
\le
\exp\left(-\sigma_0^2c_3^2M/(50b^2)\right).
$$
This shows that one column of $\mb{V}_{12}$ exceeds norm $1/2$ with
exponentially small probability.  Applying the union bound over
all the columns, we find
\begin{equation}
\mathbb{P}\left(\Vert \mb{V}_{12}\Vert_\infty\ge 1\right)
\le
(n-N)\exp\left(-\sigma_0^2c_3^2M/(50b^2)\right).
\label{eq:probv12}
\end{equation}
Thus, we have established that $\Vert \mb{V}_{12}\Vert_\infty\le 1$
with probability exponentially close to 1.

We now consider the matrix $\mb{W}_{12}$, which can be written as
$\mb{W}_{12}=\lambda^*\mb{D}\mb{R}_{12}$, where
$\mb{D}\in\R^{M\times M}$ is given by
$$\mb{D}=\mb{I}-\frac{1}{\norm{\mb{u}_1}_1}\mb{e}_M\mb{u}_1^T.$$ 
Notice that we can equivalently write
$$\mb{W}_{12}=\lambda^*\mb{D}(\mb{R}_{12}-c_3\sigma_0\mb{e}_M\mb{e}_{n-N}^T),$$
since $\mb{D}\mb{e}_M=\mb{0}$.  
The matrix $\mb{R}_{12}-c_3\sigma_0\mb{e}_M\mb{e}_{n-N}^T$ is
a subgaussian matrix scaled by $\sigma_0$.  Furthermore,
$$
\norm{\mb{D}_{12}}\leq 1+\frac{\sqrt{M}}{\norm{\mb{u}}_1}\leq 6,
$$
with high probability, since
$\Vert \mb{u}_1\Vert_1\ge M^{1/2}/5.$
Thus,
Lemma~\ref{lem:2norm}(ii) applied
to $\mb{W}_{12}/(\lambda^*\sigma_0)$, taking
$u=1/(2\lambda^*\sigma_0)$, yields
$$\mathbb{P}\left(\norm{\mb{W}_{12}}\geq 1/2\right)\leq 
\exp\left(-\left(\frac{2}{36\cdot 81b^2\sigma_0^2(\lambda^*)^2}-(\log 7)(M+n-N)\right)\right).$$
Apply the upper bound on $\lambda^*$ from \req{eq:lambdastarbounds} to 
obtain
$$\mathbb{P}\left(\norm{\mb{W}_{12}}\geq 1/2\right)\leq 
\exp\left(-\left(\frac{2(1+c_3-c_5)^2MN}{36\cdot 81b^2(1+\theta\sqrt{MN})^2}-(\log 7)(M+n-N)\right)\right).$$
Now finally we apply \req{eq:sqtheta} to obtain
\begin{equation}
\mathbb{P}\left(\norm{\mb{W}_{12}}\geq 1/2\right)\leq 
\exp\left(-\left(\frac{8c_5^2MN}{36\cdot 81b^2}-(\log 7)(M+n-N)\right)\right).
\label{eq:probw12}
\end{equation}

The same construction and analysis applies to $\mb{V}_{21}$ and
$\mb{W}_{21}$, and the same results are obtained except with the
roles of $(M,m)$ and $(N,n)$ interchanged.  
Thus,  
\begin{equation}
\mathbb{P}\left(\Vert \mb{V}_{21}\Vert_\infty\ge 1\right)
\le
(m-M)\exp\left(-\sigma_0^2c_3^2N/(50b^2)\right),
\label{eq:probv21}
\end{equation}
and
\begin{equation}
\mathbb{P}\left(\norm{\mb{W}_{21}}\geq 1/2\right)\leq 
\exp\left(-\left(\frac{8c_5^2MN}{36\cdot 81b^2}-(\log 7)(N+m-M)\right)\right).
\label{eq:probw21}
\end{equation}

{From} the analysis of all four blocks of $\mb{V}$ and $\mb{W}$, we have: 
$$
\norm{\mb{V}}_{\infty}=\max\left\{\norm{\mb{V}_{11}}_{\infty},\norm{\mb{V}_{12}}_{\infty},\norm{\mb{V}_{21}}_{\infty},\norm{\mb{V}_{22}}_{\infty}\right\}\leq 1,
$$
where $\mb{V}_{11}=\mb{e}_M\mb{e}_N^T$. With a high probability, we also have $\norm{\mb{W}}\leq 1$ since$$
\norm{\mb{W}}^2\leq\norm{\mb{W}_{11}}^2+\norm{\mb{W}_{12}}^2+\norm{\mb{W}_{22}}^2+\norm{\mb{W}_{21}}^2\le 1.
$$

By the union bound, the probability of failure of the main result
is at most the sum of the probabilities of the failure at each step.
Therefore, the failure of the convex relaxation to find the claimed
optimal $\mb{X}$ is at most the sum of the right-hand sides of 
\req{eq:probq}, \req{eq:probw22},   \req{eq:probv21},  \req{eq:probw12},
\req{eq:probv12},
and \req{eq:probw21}.  We require these probabilities
to be exponentially small.  We assure that \req{eq:probq} is 
exponentially small by requiring
\begin{equation}
MN\ge k_1(M+N)^{4/3}
\label{eq:MN1}
\end{equation}
where 
\begin{equation}
k_1> \left((\log 7)81b^2\right)^{4/3}.
\label{eq:k1def}
\end{equation}
Next, all of \req{eq:probw22}, \req{eq:probw12},
\req{eq:probw21} are exponentially small
provided that
\begin{equation}
MN\ge k_2(m+n)
\label{eq:MN2}
\end{equation}
where
\begin{equation}
k_2>\frac{(\log 7)36\cdot 81b^2}{8c_5^2}.
\label{eq:k2def}
\end{equation}
Finally, to ensure that  \req{eq:probv12} and \req{eq:probv21}
tend to 0 exponentially
fast requires that $M$ grow as fast as
$\Omega(\log(n-N))$ and similarly $N$ grows as fast as
$\Omega(\log(m-M))$, but this is already a consequence of 
\req{eq:MN1} and \req{eq:MN2}.

\section{Conclusions}
We have shown that a convex relaxation can find a large, approximately
rank-one submatrix of a much larger noisy matrix provided that the
dimensions of the larger matrix are no larger than the square of the
dimensions of the smaller matrix, and provided certain upper bounds
are satisfied on the level of the noise.

It is interesting to note that our result also applies to the
maximum biclique problem, which
was introduced in Section~\ref{sec:intro} as a special
case of LAROS.  In particular, if $G$ is a bipartite
graph $(U,V,E)$ containing a biclique given by $U^*\times V^*$, where
$|U|=m$, $|V|=n$, $|U^*|=M$, $|V^*|=N$, and if the remaining edges
of $E$ (i.e., those not in $U^*\times V^*$) are inserted at random with
probability $1/2$, then the $U$-to-$V$ adjacency matrix has the
form \req{eq:Aform} in which $\sigma =1$, $c_1=c_2=0$, $c_3=1/2$, 
$b=1/(8\log 2)^{1/2}$.  (This is not quite correct since in this
case $\mb{R}_{11}=\mb{0}$.  However, our analysis covers this case
as well.)  Thus, our algorithm with parameter $\theta=O(1/(MN)^{1/2})$
finds the planted biclique when $M\sim N$, $m\sim n$, and
$M\ge \Omega(m^{1/2})$.
The same result was obtained earlier
by Ames and Vavasis \cite{AmesVavasis} using a different convex
relaxation.  Theirs has the advantage that $M,N$ do not need to be
known or estimated in advance, but ours solves a more general class
of problems.

\bibliographystyle{plain}
\bibliography{SROFFF}
\end{document}